\documentclass[a4paper,11pt]{article}

\usepackage{a4}
\usepackage{geometry} 
\usepackage{amsmath}
\usepackage{amssymb}
\usepackage{amsthm}
\usepackage{graphicx}
\usepackage{caption,subcaption}
\usepackage{multirow}
\usepackage{xstring}
\usepackage[utf8]{inputenc}
\usepackage{hyperref}
\usepackage{amstext,amsbsy,amsopn,eucal,enumerate}
\usepackage{tikz}
\usepackage{pgfplots}
\usepackage{tabularx}

\usepackage{soul}
\usepackage{algorithm} 
\usepackage{algpseudocode}
\DeclareMathOperator*{\argmin}{arg\,min}

\def\R{\mathbb{R}}
\def\C{\mathbb{C}}

\newcommand{\expn}{\operatorname{e}}

\newcommand{\diag}{\operatorname{diag}}

\newcommand{\beq}{\begin{equation}}
\newcommand{\eeq}{\end{equation}}
\newcommand {\mat}      [1] {\left[\begin{array}{#1}}
\newcommand {\rix}          {\end{array}\right]}
\newcommand {\smat}      [1] {\left[\begin{smallmatrix}{#1}}
\newcommand {\srix}          {\end{smallmatrix}\right]}
\newcommand {\s}      [1] {\begin{smallmatrix}{#1}}
\newcommand {\se}          {\end{smallmatrix}}
\newcommand{\trace}{\operatorname{tr}}

\newtheorem{defn}{Definition}[section]
\newtheorem{remark}{Remark}

\newtheorem{lem}[defn]{Lemma}
\newtheorem{prop}[defn]{Proposition} 
\newtheorem{kor}[defn]{Corollary}
\newtheorem{thm}[defn]{Theorem}

\def\addlegendimage{\csname pgfplots@addlegendimage\endcsname}

\newlength\fheight
\newlength\fwidth

 \title{Complexity reduction of large-scale stochastic systems using linear quadratic Gaussian balancing}
 \author{Tobias Damm\thanks{RPTU Kaiserslautern-Landau, Department of Mathematics, Gottlieb-Daimler-Straße 48, 67663 Kaiserslautern, Germany, Email: {\tt 
 damm@mathematik.uni-kl.de}.}\and Martin Redmann \thanks{Corresponding author: MLU Halle-Wittenberg, Institute of Mathematics, Theodor-Lieser-Str. 5, 06120 Halle (Saale), Germany, Email: {\tt  martin.redmann@mathematik.uni-halle.de}.}
 }

\begin{document}

\maketitle

\begin{abstract}
In this paper, we consider a model reduction technique for stabilizable and detectable stochastic systems. It is based on a pair of Gramians that we analyze in terms of well-posedness. Subsequently, dominant subspaces of the stochastic systems are identified exploiting these Gramians. An associated balancing related scheme is proposed that removes unimportant information from the stochastic dynamics in order to obtain a reduced system. We show that this reduced model preserves important features like stabilizability and detectability. Additionally, a comprehensive error analysis based on eigenvalues of the Gramian pair product is conducted. This provides an a-priori criterion for the reduction quality which we illustrate in numerical experiments.
\end{abstract}



\section{Introduction}
Simulation and optimal control  of high-dimensional stochastic
processes is extremely challenging but of significant practical
interest. Such processes occur for instance as solutions to spatially discretized stochastic partial differential equations. Therefore, it is vital to reduce the computational complexity when solving such large-scale stochastic differential equation (SDE) numerically. Different techniques for model order reduction of deterministic systems have been developed over the years and are well documented,
e.g., in \cite{antoulas, BennCohe17, ObinAnde01}. A prominent method, balanced truncation, is motivated by
energy functionals and appropriate balancing of states. This is
achieved via a pair of positive definite matrices, the observability and the reachability Gramian. Given these matrices, a state space transformation can be computed such that both Gramians are equal and diagonal. Then, those
states are truncated that correspond to low output and high control
energy. Advantages of this approach are the preservation of system
properties such as stability and minimality as well as good error
bounds, as has been proved for stable systems already in \cite{moo1981,  pernebo1982model}. An extension to unstable systems using techniques from
linear quadratic (LQ) control theory has been suggested by
\cite{lqg_balancing}. Under suitable conditions the reduced unstable
system may be used in a low-order compensator to stabilize the
original system.

In the current paper, we extend this idea further to
stochastic systems. We build upon earlier work, in particular
\cite{bennerdamm, bennerdammcruz, redmannbenner},
where different versions of balanced truncation for asymptotically
stable stochastic linear systems have been discussed. The situation is
more complicated than in the deterministic setup, since frequency
domain considerations are not possible and 
hence essential tools like transfer functions are not available. Also
the duality principle of reachability and observability does not
translate literally to the stochastic setup. Therefore, it is not
immediate to find an appropriate pair of Gramians, as has been
discussed in \cite{bennerdammcruz} for the stable case.
This seems even more difficult in the unstable case, and we regard it as one of
our main contributions in this paper to suggest such a pair. While our
observability Gramian is given as the stabilizing solution of the Riccati equation
associated to a LQ-state feedback problem of stochastic
control, our reachability Gramian solves a modified Riccati type
inequality. Both Gramians exist under natural stabilizability and
detectability conditions. They can be computed, e.g., by semidefinite
programming and yield a balancing state space transformation.
Performing a balancing procedure in the usual way, we can show that the
reduced system is still stabilizable and detectable and that the truncated
closed-loop LQ-controller of the full system also stabilizes the
reduced system. Moreover, we prove error bounds for the closed-loop and
the open-loop input output system with respect to the $L^2$-norm. In
both cases, non-trivial technical elaborations are required. Further, these results can be interpreted 
nicely in the gap metric as we point out in this paper as well.

The concepts used in our approach have been developed over many years.
Fundamental results on stochastic stability can be found in
\cite{staboriginal}. Stochastic linear quadratic control theory and
the stochastic Riccati equation were introduced in \cite{wonham}. For
its solution, notions of stochastic stabilizability and detectability
are crucial. Different versions appeared in \cite{Damm07, DragHala97, LiWang10, Tess94, Tess97, WillWill76, ZhanChen04, ZhanZhan08} and have been
adapted for other classes of systems in more recent years. 
In this paper, we follow the definitions of detectability given in
\cite{Damm07, ZhanZhan08}. Low-order compensators for stochastic
linear systems with multiplicative noise apparently have first been
considered in \cite{Hinrichsen_Pritchard} using $H^\infty$-techniques.
Model order reduction of discrete-time stochastic systems based on balancing
was  discussed in \cite{ZhanHuan06} using linear matrix inequalities.
For continuous-time systems Gramian based methods were suggested in
\cite{bennerdamm, bennerdammcruz}. A balancing procedure in a Hilbert space setting can be found in \cite{beckerhartmann}. We refer, e.g., to \cite{hartmann2011} for a dimension reduction scheme based on an averaging principle. Besides  methods relying on Gramians, further recent developements in different directions have been made. In \cite{mliopt}, optimization based model reduction was studied, whereas \cite{stoch_moment_match} and \cite{pod_sde} focused on techniques based on moment matching and sampling, respectively. \smallskip

The paper is now organized as follows. We first clarify the notation and provide some tools on stochastic systems, positive operators and system theoretic
notions. In Section \ref{sec:riccati-pair}, we introduce the pair of Gramians and
characterize them with energy cost functionals. Linear quadratic Gaussian (LQG) balanced truncation is
discussed in Section \ref{sec:balanced-truncation}, where also essential preservation properties are derived. The more technical results on error bounds are given in Section \ref{sec_error} and the appendix. Some numerical
examples that illustrate and support our findings are given in Section \ref{sec_num}.

\section{Preliminaries}
\label{sec:preliminaries}

In this section, we introduce the class of stochastic systems for which we want to perform model order reduction by LQG 
balancing. To define suitable Gramians we consider the well-known Riccati equation of
the stochastic linear quadratic control problem, e.g., from
\cite{wonham}, and a new Riccati-type inequality which is inspired  by
the \emph{type II}-Gramian defined in \cite{bennerdammcruz}. We also 
recall notions of stabilizability and detectability that are essential
for the existence of the Gramians.
\subsection{Basics of stochastic systems}
We study the stochastic system
 \begin{subequations}\label{full_system}
\begin{align}\label{stochstate}
 dx(t) &= [Ax(t) +  B u(t)]\,dt + \sum_{i=1}^{q}  N_i x(t) \,dW_i(t),\quad x(0)=x_0\in\mathbb R^n,\\ \label{stoch_out}
 y(t)&= C x(t),\quad t\geq 0,
\end{align}
\end{subequations}
where $A, N_i \in \R^{n\times n}$, $B \in \R^{n\times m}$ and $C \in
\R^{p\times n}$. The vector-valued functions $x$, $u$, and $y$ are
called state, control input, and measured output respectively. We assume that
$W=\left(W_1, \ldots, W_{q}\right)^\top$ is an $\mathbb R^{q}$-valued
Wiener process with mean zero and covariance matrix $K = (k_{ij})$,
i.e., $\mathbb E[W(t)W(t)^\top] = K t$. All stochastic processes appearing in
this paper are defined on a filtered probability space 
$\left(\Omega, \mathcal F, (\mathcal F_t)_{t\geq 0}, \mathbb
  P\right)$\footnote{$(\mathcal F_t)_{t\geq 0}$ is right continuous
  and complete.}. Furthermore, we assume that $W$ is an $(\mathcal F_t)_{t\geq
  0}$-adapted process with  
increments $W(t+h)-W(t)$ being independent of $\mathcal F_t$ for $t,
h\geq 0$. Throughout this paper, suppose that $u$ is an $(\mathcal
F_t)_{t\geq 0}$-adapted control with $u\in L^2_T$, meaning that  
\begin{align}\label{control_assumption}
\left\| u\right\|^2_{L^2_T} := \mathbb E\int_0^T \left\| u(s)\right\|_2^2 ds < \infty
\end{align}
for all $T>0$, where $\left\| \cdot\right\|_2$ denotes the Euclidean norm with associated inner product $\langle \cdot, \cdot\rangle_2$. If \eqref{control_assumption} additionally holds for $T=\infty$,
we write $u\in L^2$. For given control $u$ and initial state
$x_0$, the corresponding state and output processes are denoted by
$x(\cdot, x_0, u)$ and $y(\cdot, x_0, u)$.
\begin{remark}
  There is the potential to extend the results of this paper to square
  integrable L\'evy processes (see, e.g., \cite{redmannPhD}). We might   also consider multiplicative noise at the input terms. This, however, makes many expressions and criteria more complicated, as we demonstrate below in Remark \ref{tobeadded}.
\end{remark}

\subsection{Resolvent positive mappings }
In our analysis we will consider Lyapunov equations of a generalized
type. In this context the following terminology and setup is useful, see
\cite{damm}.
Let $H$ denote a finite dimensional real vector space ordered by 
a closed, solid, pointed convex cone $H_+$.
A linear mapping $T:H\to H$ is called positive,
if $T(H_+)\subset H_+$. It is called
\emph{resolvent positive}, if its resolvent $(\alpha I-T)^{-1}$ is
positive for all sufficiently large $\alpha$. The essential property
of resolvent positive mappings that we use is a variant of the
Perron-Frobenius theorem.
\begin{prop}\label{prop:PF}
  Let $T:H\to H$ be resolvent positive with
  spectrum $\sigma(T)$ and spectral abscissa
  $\alpha=\max\{\Re(\lambda)\;\big|\;\lambda\in\sigma(T)\}$.
  Then $\alpha\in\sigma(T)$ and there exists $X\in H_+$,
  $X\not=0$, such that $T(X)=\alpha X$.
\end{prop}
In our context, we consider the space $\mathcal{S}^n=\{X\in\mathbb{R}^{n\times n}\;\big|\;
X=X^\top\}$ of symmetric matrices. This space is endowed with the canonical Frobenius scalar product  $\langle
X_1,X_2\rangle_F=\trace(X_1X_2)$ and ordered by the closed, solid, pointed convex cone of nonnegative definite
matrices $\mathcal{S}^n_+=\{X\in\mathcal{S}^n\;\big|\; X\ge0\}$.  We will use the following property.
\begin{align}
  \label{eq:sp0}
 \text{If } X_1,X_2\in \mathcal{S}^n_+ \text{ then } X_1X_2=0 \iff
  \langle X_1,X_2\rangle_F =0\;.
\end{align}

With the given coefficient matrices $A,N_i\in\mathbb{R}^{n\times n}$, 
$K=(k_{ij})\in\mathcal{S}^q_+$ from the previous subsection we define the mappings
$\mathcal{L}_A,\Pi_N:\mathcal{S}\to \mathcal{S}$ by
\begin{align*}
  \mathcal{L}_A(X)=A^\top X+XA\;,\quad \Pi_N(X)=\sum_{i,j=1}^q
  N_i^\top X N_j\;k_{ij}.
\end{align*}
Then, $\Pi_N$ is positive and the sum $\mathcal{L}_A+\Pi_N$ is
\emph{resolvent positive}. The same obviously holds for the adjoint mappings
\begin{align*}
  \mathcal{L}_A^*(X)=A X+XA^\top\;,\quad \Pi_N^*(X)=\sum_{i,j=1}^q
  N_i X N_j^\top\;k_{ij}.
\end{align*}

\subsection{Stabilizability, observability and detectability}
\label{sec:stab-detect}

We now introduce notions of stability, stabilizability, and
detectability, as they have been considered, e.g.,  in \cite{damm,
  Damm07}. 
\begin{defn}\label{def:stab_obs_dtct}
  The system \eqref{full_system} is called
  \begin{itemize}
  \item \emph{mean square asymptotically stable}, if there exist
    constants $M\geq 1, c>0$ such that for all $x_0\in \mathbb R^n,
    t\ge 0$, we have $\mathbb E\left\|x(t, x_0, 0)\right\|_2^2 \leq M \expn^{-c t}\left\|x_0\right\|_2^2$.
  \item \emph{stabilizable}, if for all   $x_0\in \mathbb R^n$ there
    exists $u\in L^2$, such that $x(\cdot, x_0, u)\in L^2$.
    \item \emph{observable}, if the condition that $y(t, x_0, 0)=0$ almost
     surely for all $t\ge 0$ implies that $x_0=0$. 
  \item \emph{detectable}, if the condition that $y(t, x_0, 0)=0$ almost
    surely for all $t\ge 0$ implies that $\lim_{t
      \to\infty}\mathbb E\left\|x(t, x_0, 0)\right\|_2^2=0$.
  \end{itemize}
 In these cases, we briefly say that the pair $(A,N_i)$ is (mean square asymptotically) stable, the
 triple $(A,B,N_i)$
 is stabilizable, or the triple $(A,C,N_i)$ is observable or detectable.
\end{defn}

\begin{remark}
  In \cite{ZhanChen04} the term \emph{exact observable} is used,
  where for brevity we just write \emph{observable}. To be more precise, the notion of stability introduced in Definition \eqref{def:stab_obs_dtct} is called mean square exponential stability in general. Since usual mean square asymptotic stability implies exponentially fast decay in the linear case, we do not distiguish between both concepts and omit the term ``exponential'' in the following.
\end{remark}
Unlike in the deterministic case there is no perfect duality between
stabilizability and detectability.  The following lemma collects known
criteria. 
\begin{lem}\label{properties_stoch_sys}
  \begin{itemize}
  \item[(a)] The triple $(A,B,N_i)$ is stabilizable, if and only if there
    exists a feedback gain
    matrix $F$, such that $(A+BF,N_i)$ is stable.
    \item[(b)]    The triple $(A,C,N_i)$ is observable, if and only if the following
    variant of the Hautus test is satisfied: 
    \begin{align*}
      \text{If } 
      (\mathcal{L}_A+\Pi_N)^*(V)=\lambda V\text{ with }
      \lambda \in\mathbb C , 0\neq V\ge0, \text{ then } CV\neq 0\;. 
    \end{align*}
 \item[(c)]    The triple $(A,C,N_i)$ is detectable, if and only if the following
    variant of the Hautus test is satisfied:
    \begin{align*}
      \text{If } 
      (\mathcal{L}_A+\Pi_N)^*(V)=\lambda V\text{ with }
      \lambda\ge 0, 0\neq V\ge0, \text{ then } CV\neq 0\;. 
    \end{align*}
    \item[(d)] If the triple $(A^\top,C^\top,N_i^\top)$ is
      stabilizable, then the triple $(A,C,N_i)$ is detectable. The
      converse does not hold in general.
  \end{itemize}
\end{lem}

\section{A pair of Gramians}
\label{sec:riccati-pair}

As in the deterministic case, stabilizability, observability and detectability
characterize the solvability of Riccati equations.

\subsection{An observability Gramian}
\label{sec:an-observ-gram}
We first consider the Riccati equation of the stochastic linear
quadratic control problem. The following result is a special case of \cite[Theorem 4.1]{ZhanZhan08} (see also \cite[Corollary
  5.3.4]{damm}).
\begin{thm}\label{thm:StochLQ}
  Assume that $(A,B,N_i)$ is stabilizable and  $(A,C,N_i)$ is
  detectable. Then, the Riccati equation
  \begin{align}
  \mathcal{R}(Q)&:=A^\top Q + Q A + \sum_{i, j=1}^q N_i^\top Q N_j k_{ij}+ C^\top C -QB
    B^\top Q = 0 \label{obs_gram}
  \end{align}
  possesses a stabilizing solution $Q_+\ge 0$, such that
  $(A-BB^\top Q_+,N_i)$ is stable.\\
  If $(A, C, N_i)$ is observable, then $Q_+>0$.
\end{thm}
The stabilizing solution $Q_+$ of \eqref{obs_gram} will play the role
of an observability Gramian in our LQG balanced truncation approach.

\subsection{A reachability Gramian}
\label{sec:reachability-gramian}

The corresponding reachability Gramian will be chosen as
a positive definite solution of the new Riccati-type inequality
\begin{align}\label{reach_gram}
 A^\top P^{-1} + P^{-1} A + \sum_{i, j=1}^q N_i^\top P^{-1} N_j k_{ij}
  -C^\top C +P^{-1}B B^\top P^{-1}\leq 0\;.
\end{align}
\begin{lem}
  The \emph{strict} inequality
  \begin{align}\label{reach_gram_strict}
 A^\top P^{-1} + P^{-1} A + \sum_{i, j=1}^q N_i^\top P^{-1} N_j k_{ij}
  -C^\top C +P^{-1}B B^\top P^{-1}< 0
  \end{align}
  possesses a solution $P_+>0$ if and only if
  $(A^\top,C^\top,N_i^\top)$ is  stabilizable.
\end{lem}
\begin{proof}
 By \cite[Lemma 1.7.3]{damm} stabilizability of $(A^\top, C^\top,
 N_i^\top)$ is equivalent to the existence of a matrix $X>0$ such
 that
 \begin{align}\label{ineq:Lemma173}
A^\top X + X A + \sum_{i, j=1}^q N_i^\top X N_i k_{ij}- C^\top
 C = -Y<0\;.
 \end{align}
Since \eqref{reach_gram_strict}  implies \eqref{ineq:Lemma173} it also
implies stabilizability of $(A^\top, C^\top,
 N_i^\top)$.\\
For the converse implication, we multiply both sides of
\eqref{ineq:Lemma173} by some $0<
 \epsilon \leq 1$ and obtain
 \begin{align*}
   A^\top (\epsilon X) + (\epsilon X) A + \sum_{i, j=1}^q N_i^\top
   (\epsilon X) N_i k_{ij}- \epsilon C^\top C = -\epsilon
   Y.                                                                                                                                                                                                                            \end{align*}
 The choice of $\epsilon\le 1$ yields $- C^\top C\leq - \epsilon C^\top C$.
 Moreover, for sufficiently small $\epsilon$, we have
 $-\epsilon Y < -(\epsilon X) B B^\top (\epsilon X)$. Consequently,
 the corresponding $P_+ := (\epsilon X)^{-1}$ is a positive definite
 solution to  \eqref{reach_gram_strict}. 
\end{proof}
\begin{remark}
  As noted before, stabilizability of  $(A^\top, C^\top,
  N_i^\top)$ is stronger than the more natural detectability of $(A,
  C,N_i)$ and it is also not implied by observability of $(A,
  C,N_i)$, \cite{Damm07}.  However, in the extreme case, where
  $BB^\top $ is nonsingular, it is clear that \eqref{reach_gram} implies 
  \eqref{ineq:Lemma173} and thus stabilizability of  $(A^\top, C^\top,
  N_i^\top)$. In the following, we will make the assumptions that
  $(A,B,N_i)$ is stabilizable, $(A, C,N_i)$ is observable, and that \eqref{reach_gram}
  has a solution $P>0$.  
\end{remark}

\subsection{State cost estimations}
\label{sec:state-cost-estim}

In this section, we measure how much state variables contribute to system \eqref{full_system} based on the proposed Gramians.
\paragraph{Closed-loop dynamics}
The relevance of state components with respect to
 the quadratic cost functional
 \begin{align*}
   J_T(x_0,u) =\int_0^T \mathbb{E}\left(\|u(t)\|_2^2+\|y(t)\|_2^2\right)\,dt
 \end{align*}
is investigated. Let us first assume that $x_0=0$ and hence $x(t)=x(t, 0, u)$. Moreover, suppose that we have an orthonormal basis $(p_i)$ of eigenvectors $P>0$ such that we have the representation 
\begin{align}\label{state_expansion}
 x(t) = \sum_{i=1}^n \left\langle x(t),  p_i \right\rangle_2 p_i.
\end{align}
In order to tell how much a direction $p_i$ contribute to the state variable, the coefficients $\left\langle x(t),  p_i \right\rangle_2$ are analyzed below. Secondly, we investigate how much a state variable contributes to the output (and a feedback control). Since a state is fully determined by its initial condition, we focus on $x_0= \sum_{i=1}^n \beta_i q_i$. Here, $(q_i)$ is an orthonormal basis of eigenvectors of $Q$ and $\beta_i\in\mathbb R$ and the respective coefficients of the expansion of $x_0$. We assume that the control has stabilizing feedback structure, i.e.,  $u(t, x_0)= - B^\top Q x(t, x_0)$. Since the state is linear in $x_0$ the same is true for the output $y=y(t, x_0)$ and the closed-loop control. Consequently, we have\begin{align}\label{output_expansion}
y(t, x_0)=  \sum_{i=1}^n \beta_i y(t, q_i)\quad\text{and}\quad  u(t, x_0)=  \sum_{i=1}^n \beta_i u(t, q_i).                                                                                                                                                                                                                                                       \end{align}
Therefore, it is of interest to investigate how large $u(t, q_i)$ and $y(t, q_i)$ are. We establish the following proposition in order to characterize dominant subspaces.
\begin{prop}\label{energy_est}
   \begin{itemize}
   \item [(a)] If $P>0$ satisfies
     \eqref{reach_gram}, and $x(t)=x(t,0,u)$ is the solution of
     \eqref{stochstate} on $[0,T]$ with initial value $x(0)=0$, then we have
\begin{align*}
  \sup_{t\in[0,T]}\mathbb{E}\left\langle x(t),  p_i \right\rangle_2^2 &\le \lambda_{P,i} J_T(0,u)
\end{align*}
for the coefficients in \eqref{state_expansion}, where $\lambda_{P,i}$ is the eigenvalue of $P$ associated to $p_i$. 
If $u$ is stabilizing, then the same holds for $T=\infty$.
\item[(b)] Let now $Q\geq 0$ satisfy
     \eqref{obs_gram}, and consider the initial state $x_0=q_i$
     with stabilizing feedback input $u_F=Fx$ for $F=-B^\top Q$.
    Then, \begin{align}\label{inequality_observability}                                                                                                                                J_\infty(q_i,u_F) = \left\| y(\cdot, q_i)\right\|^2_{L^2}+\left\| u_F(\cdot, q_i)\right\|^2_{L^2} \leq \lambda_{Q, i}                                                                                                                            \end{align}
for the coefficients in \eqref{output_expansion}, where $\lambda_{Q, i}$ is the eigenvalue of $Q$ associated to $q_i$.
\end{itemize}
 \end{prop}
 \begin{remark}\label{remark_dom}
 Note that each $u\in L^2$ is stabilizing in Proposition \ref{energy_est} (a) if \eqref{stochstate} is mean square asymptotically stable \cite{Hinrichsen_Pritchard}.
   We interpret the results of the proposition as follows. Assume that $\lambda_{P, i}$ and $\lambda_{Q, i}$ are very small. Then, from (a) we infer, that the state direction
    $p_i$ can only be activated at very high cost. From (b)
   we learn that the initial variable direction $q_i$ has only very small
   influence on the cost meaning that $u_F(\cdot, q_i)$ and $y(\cdot, q_i)$ are small in $L^2$. If $P$ and $Q$ are diagonal and equal, we obtain that $p_i=q_i$ are the unit vectors such that unimportant directions can be identified with state components. The $i$-th component might then be neglected in a truncation approach if the associated diagonal entry in $P=Q$ is small. Such a simultaneous diagonalization of the Gramians
   is discussed in Section \ref{bal_trans}.
 \end{remark}
\begin{proof}[Proof of Proposition \ref{energy_est}]
For proving (a), we use \eqref{estimate_quadr_functional} with $X=P^{-1}$ and \eqref{reach_gram} yielding 
\begin{align*}
 &\mathbb E \langle x(t, 0, u), p_{i}  \rangle_2^2 \leq 
  \lambda_{P, i} \;\mathbb E \left[x(t, 0, u)^\top P^{-1} x(t, 0, u)\right]\\
  &\leq  \lambda_{P, i} \left[\int_0^t \mathbb E\left[x(s)^\top\left(C^\top C- P^{-1} B B^\top P^{-1} \right) x(s)\right] ds+ 2 \int_0^t \mathbb E\left\langle B^\top P^{-1} x(s),  u(s) \right\rangle_2ds\right]\\
&=\lambda_{P, i} \left[\left\| y\right\|^2_{L^2_t}+\left\| u\right\|^2_{L^2_t} - \left\| B^\top P^{-1} x-u\right\|^2_{L^2_t}\right].
\end{align*}
Consequently, we obtain
\begin{align*}
 \sup_{t\in[0, T]} \mathbb E \langle x(t, 0, u), p_{i}  \rangle_2^2 \leq \lambda_{P, i}\left[\left\| y\right\|^2_{L^2_T}+\left\| u\right\|^2_{L^2_T}\right].
\end{align*}
If $u, y\in L^2$, we can take the supremum over $t\in [0, \infty)$ instead. For (b), we set $X=Q$ in \eqref{estimate_quadr_functional} and make use of \eqref{obs_gram} leading to
\begin{align}\nonumber
\mathbb E \left[x(t)^\top Q x(t)\right] &= 
x_0^\top Q x_0+\int_0^t \mathbb E\left[x(s)^\top\left(- C^\top C +QB B^\top Q \right) x(s)\right] ds\\ \nonumber
&\quad + 2 \int_0^t \mathbb E\left\langle B^\top Q x(s),  u(s) \right\rangle_2ds\\ \label{estimate_on_v}
&=x_0^\top Q x_0+\int_0^t \mathbb E\left[-\left\|y(s)\right\|_2^2-\left\|u(s)\right\|_2^2 + \left\| B^\top Q x(s) + u(s)\right\|_2^2\right] ds
\end{align}
for $t\in[0, T]$. 
Setting $u=u_F$ and $x_0=q_i$ gives us \begin{align*}                                                                                                                                \left\| y(\cdot, q_i)\right\|^2_{L^2_T}+\left\| u_F(\cdot, q_i)\right\|^2_{L^2_T} \leq \lambda_{Q, i}.                                                                                                                           \end{align*}
Since $u_F$ is a stabilizing control according to Theorem \ref{thm:StochLQ}, the result follows by taking $T\rightarrow \infty$.
\end{proof}
\begin{remark}\label{tobeadded}
  Using the Riccati mapping $\mathcal{R}$ from \eqref{obs_gram}, we
  can write \eqref{reach_gram} in the form $-\mathcal{R}(-P^{-1})\le
  0$. This indicates, how the Gramians can be generalized for the case
 of models with control-dependent noise. Given $M_i\in \mathbb R^{n\times m}$, let us
  consider the following system with controlled diffusion
\begin{align*}
 dx(t) &= [Ax(t) +  B u(t)]\,dt + \sum_{i=1}^q[N_i x(t) +M_i u(t)] \,dW_i(t),\quad x(0)=x_0,\\ 
 y(t)&= C x(t),\quad t\geq 0.
\end{align*}
Then, the Riccati equation of LQ-control with cost functional
$J_\infty$ takes the form (see, e.g., \cite{Bism76, damm})
\begin{align*}
  \mathcal{R}(Q)&=A^\top Q+QA+\sum_{i, j=1}^q N_i^\top QN_j k_{ij}+C^\top C\\&-(QB+\sum_{i, j=1}^q N_i^\top Q M_j k_{ij})(I+\sum_{i, j=1}^q  M_i^\top Q M_j k_{ij})^{-1}(B^\top Q+\sum_{i, j=1}^q  M_i^\top Q N_j k_{ij})=0,
\end{align*}
and defines the observability Gramian for this case. A reachability
Gramian is given by
\begin{align*}
  0&\geq -\mathcal{R}(-P^{-1})\\
  &=A^\top P^{-1}+P^{-1}A+\sum_{i, j=1}^q N_i^\top P^{-1}N_j k_{ij}-C^\top C\\&+(P^{-1}B+\sum_{i, j=1}^q N_i^\top P^{-1}M_jk_{ij})(I-\sum_{i, j=1}^q M_i^\top P^{-1}
  M_j k_{ij})^{-1}(B^\top P^{-1}+\sum_{i, j=1}^q M_i^\top P^{-1}N_j k_{ij}),
\end{align*}
if additionally 
\begin{align}\label{add_constr}
I-\sum_{i, j=1}^q M_i^\top P^{-1}M_j k_{ij} > 0.
\end{align}
Proposition \ref{energy_est}
holds accodingly in this setup with $F=-(I+\sum_{i, j=1}^q  M_i^\top Q M_j k_{ij})^{-1}(B^\top Q+\sum_{i, j=1}^q  M_i^\top Q N_j k_{ij})$, but the additional constraint \eqref{add_constr} on $P$ causes further technical difficulties. Therefore, we prefer to consider
only state-dependent noise.
\end{remark}
\paragraph{Open-loop dynamics}
Below, we discuss that the Gramian $Q$ might not generally be suitable for the dominant subspace characterization of unstable (but stabilizable and detectable) systems. Let $u$ now be an open-loop control. Since $u$ then is independent of the (initial) state, we can neglect it in the considerations below by setting $u\equiv 0$, i.e., $y(t) = y(t, x_0, 0) = C x(t, x_0, 0)$. Equation \eqref{estimate_on_v} yields \begin{align*}
\mathbb E \left[x(t)^\top Q x(t)\right] =a(t) + \int_0^t\left\| B^\top Q x(s) \right\|_2^2 ds,\quad t\in[0, T],
\end{align*}
where $a(t):=x_0^\top Q x_0-\left\|y\right\|_{L_t^2}^2$. Using $\left\| B^\top Q x(s) \right\|_2^2\leq b\;x(s)^\top Q x(s)$ with $b:= \left\| B^\top Q^{\frac{1}{2}} \right\|_2^2$, we can apply Gronwall's lemma, see Lemma \ref{gronwall}. Setting $t=T$ we then obtain \begin{align}\nonumber
\mathbb E \left[x(T)^\top Q x(T)\right] &=a(T) + \int_0^T a(s) b \expn^{b(T-s)}ds= a(T) - \left[ a(s)\expn^{b(T-s)}\right]_{s=0}^T + \int_0^T \dot a(s) \expn^{b(T-s)}ds\\ \label{open_loop_estimate}
&=x_0^\top Q x_0\expn^{bT} - \int_0^T \left\| y(s) \right\|_2^2 \expn^{b(T-s)}ds.
\end{align}
We have  $y(t, x_0, 0)=  \sum_{i=1}^n \beta_i y(t, q_i, 0)$, $t\in [0, T]$. We obtain from \eqref{open_loop_estimate} that the $i$-th summand of this expansion satisfies 
 \begin{align}\label{obs_open_loop}
\int_0^T \expn^{-bs} \left\| y(t, q_i, 0) \right\|_2^2 ds\leq \lambda_{Q, i}.
\end{align}
Inequality \eqref{obs_open_loop} is based on a Gronwall estimate that generally is not tight but captures the worst-case scenarios. Therefore, 
the exponential weight in \eqref{obs_open_loop} is an indicator that eigenspaces corresponding to small eigenvalues of $Q$ might generally only be redundant in an unstable open-loop system on a small time scale. However, the kernel of $Q$ remains negligible in any case. We proceed with a strategy to simultaneous diagonalize $P$ and $Q$ in order to be able to remove redundant information in \eqref{stochstate} and \eqref{stoch_out} at the same time.
 \subsection{State space transformation and balancing}\label{bal_trans}
A transformation $z=Sx$ in \eqref{full_system} with nonsingular
$S\in\mathbb{R}^{n\times n}$ leads to a an equivalent stochastic
system with the state vector $z$, where the coefficient matrices
undergo the state-space transformation
\begin{align*}
  (A,N_i,B,C)\mapsto (SAS^{-1},SN_iS^{-1},SB,CS^{-1}).
\end{align*}
Both systems have the same input and output. Also, none of the
properties from Definition \ref{def:stab_obs_dtct} is affected.
Matrices $P$ and $Q$ constitute a pair of Gramians for the original system, if
and only if $SPS^\top$ and $S^{-\top}QS^{-1}$ constitute a
pair of Gramians for the transformed system.
By the spectral transformation theorem for symmetric matrices, there
exist orthogonal matrices $S_P$ and $S_Q$, such that
$S_PPS_P^\top=\Sigma_P$ and $S^{-\top}_QQS^{-1}_Q=\Sigma_Q$ are
diagonal and contain the ordered eigenvalues of $P$ and $Q$, respectively.

Given that $P, Q>0$ it is possible to conduct a balancing procedure, where one computes a
nonsingular (but not necessarily orthogonal)
transformation matrix $S_b$, so that $S_bPS_b^\top=S_b^{-\top}QS_b^{-1}=\Sigma_n$
is diagonal, with $\Sigma^2_n=\diag(\sigma_1^2,\ldots,\sigma_n^2)>0$
containing the ordered eigenvalues of $PQ$. One way
of choosing the matrix $S_b$ is to compute a Cholesky factorization of
 $P=L_PL_P^\top$ and then a spectral factorization of $L_P^\top QL_P=U\Sigma_n^2 U^\top$
 with orthogonal $U$, where $\Sigma_n$ turns out to be the balanced
 Gramian. According to Remark \ref{remark_dom}, state components associated to small diagonal entries of $\Sigma_n$ are less relevant in a balanced system. They can be removed due to their low contribution to the dynamics. This idea is the basis for the reduced model introduced in the next section.

 \section{LQG balanced truncation}
\label{sec:balanced-truncation}
Our standing assumption is that 
  $(A,B,N_i)$ is stabilizable, $(A, C,N_i)$ is observable, and that \eqref{reach_gram}
  has a solution $P>0$.  
  Then, also \eqref{obs_gram} has a stabilizing solution $Q_+>0$. In this case, we can apply the balancing transformation $S_b$, leading to the 
  balanced realization $(A_n,N_{i, n}, B_n, C_n)= (S_bAS_b^{-1},S_bN_iS_b^{-1},S_bB,CS_b^{-1})$ with diagonal Gramians
  $S_b PS_b^\top=S_b^{-\top}Q_+S_b^{-1}=\Sigma_n=\diag(\Sigma_r,\Sigma_{2, n-r})$, where
  $\Sigma_r= \diag(\sigma_1,\ldots,\sigma_r)$ contains the large and $\Sigma_{2, n-r}=\diag(\sigma_{r+1},\ldots,\sigma_n)$, $r<n$, the small singular values. The balanced system matrices are partitioned conformingly
  \begin{align}\label{part_bal}
A_n= \begin{bmatrix}{A}_{r}&\star\\ 
\star&\star\end{bmatrix},\quad B_n = \begin{bmatrix}{B}_r\\\star\end{bmatrix}, \quad C_n= \begin{bmatrix}{C}_r &
\star\end{bmatrix},\quad N_{i, n}= \begin{bmatrix}{N}_{i, r}&\star\\ 
\star&\star\end{bmatrix}\;.
  \end{align}
 Then, we consider the reduced system
\begin{subequations}\label{red_system}
\begin{align}\label{red_stochstate_mul}
 dx_r(t) &= \left[A_r x_r(t) +  {B_r} u(t)\right]\,dt + \sum_{i=1}^{q} {N}_{i, r} x_r(t) \,dW_i(t),\quad x_r(0) =  x_{0, r}\in\mathbb R^r,\\
 y_r(t) &= C_r x_r(t),\quad t\geq 0.
\end{align}
\end{subequations}
By Theorem \ref{thm:StochLQ} the closed-loop system $(A_n-B_n B_n^\top\Sigma_n,N_{i, n})$  is
stable. We will show that the same holds for the reduced closed-loop
system $(A_r-B_r B_r^\top\Sigma_r,N_{i,r})$, if $\sigma(\Sigma_r)\cap\sigma(\Sigma_{2, n-r})=\emptyset$.  
Moreover, we prove that also detectability
is preserved by truncation.
\subsection{Preservation of closed-loop stability}
We make use of a result in
\cite{bennerdammcruz, dammbennernewansatz}, which we restate here in a
suitable form, see \cite[Theorem II.2]{ bennerdammcruz}.
\begin{thm}\label{thm:dualGrams_restated}
  Let $(\hat A,\hat N_i,\hat B,\hat C)$ be coefficient matrices with
  the same partitioning as in \eqref{part_bal}.
  Assume that $(\hat A,\hat N_i)$ is stable and
consider the systems
\begin{subequations}\label{eq:full_reduced_hat}
  \begin{align}
    d \hat x(t) &= \left[\hat A \hat x(t) +  \hat B \hat u(t)\right]dt + \sum_{i=1}^{q}  \hat N_i \hat x(t)
          dW_i(t),\quad
         \hat y(t)= \hat  C  \hat x(t),\\ \label{rom_bt}
    d \hat x_r(t) &= \left[\hat A_r \hat x_r(t) +  \hat B_r \hat u(t)\right]dt + \sum_{i=1}^{q}  \hat N_{i,r}
            \hat x_r(t) dW_i(t),\quad
            \hat y_r(t)= \hat C_r  \hat x_r(t).
  \end{align}
\end{subequations}
Let further $\hat\Sigma=\diag(\hat \Sigma_r,\hat\Sigma_{2, n-r})$ with
  $\sigma(\hat{\Sigma}_r)\cap\sigma(\hat{\Sigma}_{2, n-r})=\emptyset$
  satisfy 
  \begin{align}\label{eq:typeIIgrams}
    \left(\mathcal{L}_{\hat A}+\Pi_{\hat N}\right)(\hat\Sigma)&\le
                                        -\hat C^\top\hat
                                                                C\text{
                                                         and }
                                                         \left(\mathcal{L}_{\hat
                                                         A}+\Pi_{\hat
                                                                N}\right)(\hat\Sigma^{-1})\le
                                                                -\hat\Sigma^{-1}\hat
                                                                B\hat
                                                                B^\top\hat
                                                                \Sigma^{-1}\;.
  \end{align}
  Then, $(\hat A_r,\hat N_{i,r})$ is stable.
\end{thm}
\begin{remark}
  Matrices $\hat \Sigma$ satisfying \eqref{eq:typeIIgrams} have been called
  \emph{type-II-Gramians} of system $(\hat A,\hat N_i,\hat B,\hat
  C)$ and \eqref{rom_bt} the reduced model by type-II balancing, see \cite{bennerdammcruz}. 
\end{remark}
\begin{thm}\label{thm:red_stable_bound}
Consider the systems \eqref{full_system} and \eqref{red_system}  given
by the data \eqref{part_bal}. Then, we have that 
  $( A_r-B_rB_r^\top \Sigma_r,N_{i,r})$ is stable.
\end{thm}
\begin{proof}
We balance the system in order to work with the coefficient in \eqref{part_bal}. We add \eqref{obs_gram} and \eqref{reach_gram} with $P=Q=\Sigma_n$ to obtain
  \begin{align*}
    0&\ge
       \mathcal L_{A_n} (\Sigma_n+\Sigma_n^{-1})+\Pi_{N_n}(\Sigma_n+\Sigma_n^{-1})+\Sigma_n^{-1}B_nB_n^\top \Sigma_n^{-1}-\Sigma_n
       B_nB_n^\top \Sigma_n\\
     &\ge\mathcal L_{A_n-B_nB_n^\top \Sigma_n}(\Sigma_n+\Sigma_n^{-1})+\Pi_{N_n}(\Sigma_n+\Sigma_n^{-1})+\Sigma_n^{-1}B_nB_n^\top \Sigma_n^{-1}-\Sigma_n
       B_nB_n^\top \Sigma_n\\
     &\quad +\Sigma_n B_nB_n^\top  \Sigma_n+\Sigma_n B_nB_n^\top  \Sigma_n^{-1}+\Sigma_n^{-1} B_nB_n^\top  \Sigma_n +\Sigma_n B_nB_n^\top  \Sigma_n\\
     & = (\mathcal{L}_{A_n-B_nB_n^\top \Sigma}+\Pi_{N_n})(\Sigma_n+\Sigma_n^{-1})    +(\Sigma_n+\Sigma_n^{-1})B_nB_n^\top (\Sigma_n+\Sigma_n^{-1})\;.
  \end{align*}
  Let us set $\Upsilon_n=(\Sigma_n+\Sigma_n^{-1})^{-1}$. Then, we have the two inequalities
  \begin{align*}
    (\mathcal{L}_{A_n-B_nB_n^\top\Sigma_n}+\Pi_{N_n})(\Upsilon_n^{-1})&\le -\Upsilon_n^{-1} B_nB_n^\top\Upsilon_n^{-1},\\
    (\mathcal{L}_{A_n-B_nB_n^\top\Sigma_n}+\Pi_{N_n})(\Sigma_n)&\le -C_n^\top C_n-\Sigma_n B_nB_n^\top\Sigma_n.
  \end{align*}
   We recognize $\Upsilon_n$ and $\Sigma_n$ as unbalanced type-II Gramians of the closed-loop
  system given by $\left(A_n-B_nB_n^\top \Sigma_n,N_{i, n},B_n,\left[
    \begin{smallmatrix}
     -B_n^\top\Sigma_n\\ C_n
    \end{smallmatrix}
\right]\right)$. These are balanced by the similarity transformation with
  \begin{align*}
  S_n=(\Upsilon_n^{-1}\Sigma_n)^{1/4}=\diag(S_r,S_{2, n-r})\;.
  \end{align*}
  For the given $\Upsilon_n$, the balanced type-II Gramian of the closed-loop system
  then equals
  \begin{align}\label{eq:Sigma_hat}
   \hat \Sigma_n= (I+\Sigma_n^{-2})^{-1/2}=\diag\left(\frac{\sigma_j}{(1+\sigma_j^2)^{1/2}}\right)
    _{j=1}^n=\diag(\hat{\sigma}_j)_{j=1}^n\;.  
  \end{align}
  Note that $\hat\sigma_j>\hat\sigma_k$,
  if and only if $\sigma_j>\sigma_k$. Hence
  $\sigma(\Sigma_r)\cap\sigma(\Sigma_{2, n-r})=\emptyset$ implies 
  $\sigma(\hat\Sigma_r)\cap\sigma(\hat\Sigma_{2, n-r})=\emptyset$.
  Thus, the assumptions of Theorem \ref{thm:dualGrams_restated} are
  satisfied with
  \begin{subequations}\label{eq:hat_coeffs}
    \begin{align}
      (\hat A,\hat N_i,\hat B,\hat
      C)&=(S_n(A_n-B_nB_n^\top\Sigma_n)S_n^{-1},S_n N_{i, n}S_n^{-1},S_nB_n, \left[
    \begin{smallmatrix}
      -B_n^\top\Sigma\\C_n
    \end{smallmatrix}
\right] S_n^{-1}),\\
      (\hat A_r,\hat N_{i,r},\hat B_r,\hat
      C_r)&=(S_r(A_r-B_rB_r^\top\Sigma_r)S_r^{-1},S_rN_{i,r}S_r^{-1},S_rB_r, \left[
    \begin{smallmatrix}
      -B_r^\top\Sigma_r\\C_r
    \end{smallmatrix}
\right] S_r^{-1})\;.
    \end{align}
  \end{subequations}
 The stability of $(A_r-B_rB_r^\top \Sigma_r,N_{i,r})$ now follows 
 from Theorem \ref{thm:dualGrams_restated}.
\end{proof}

\subsection{Preservation of detectability}

Let us now show that the reduced system is also detectable.
\begin{prop}
  If $\sigma(\Sigma_r)\cap\sigma(\Sigma_{2, n-r})=\emptyset$, then
  $(A_r, C_r, N_{i,r})$ given by \eqref{part_bal} is detectable.
\end{prop}
\label{sec:pres-detect}
\begin{proof}
  Let us consider the balanced realization with partition in \eqref{part_bal}, so that $P=Q=\Sigma_n$ in \eqref{obs_gram} and \eqref{reach_gram}. In more detail, we partition $N_{i, n}=\left[
    \begin{array}{cc}
      N_{i,r}&\star\\M_{i,r}&\star
    \end{array}
\right]$ and define $\Pi_{M_r}:\mathcal{S}^{n-r}\to\mathcal{S}^r$ in
analogy to $\Pi_N$ by $\Pi_{M_r}(X)=\sum_{i,j=1}^q M_{i,r}^\top XM_{j,r}k_{ij}$.   Then,
  \begin{align}\label{eq:Pgram_trunc}
    \left(\mathcal{L}_{A_r}+\Pi_{N_{r}}\right)\left(\Sigma_r^{-1}\right)
   &\le
     C_r^\top C_r-\Sigma_r^{-1}B_rB_r^\top \Sigma_r^{-1}-\Pi_{M_{r}}\left(\Sigma_{2, n-r}^{-1}\right),\\\label{eq:Qgram_trunc}
       \left(\mathcal{L}_{A_r}+\Pi_{N_{r}}\right)\left(\Sigma_r\right)
   &= -C_r^\top C_r+\Sigma_rB_rB_r^\top \Sigma_r-\Pi_{M_{r}}\left(\Sigma_{2, n-r}\right).
  \end{align}
Recall that $\langle \cdot, \cdot\rangle_F$ is the Frobenius inner product. Assume that $(A_r,N_{i,r},C_r)$ is not detectable. Then, according to Lemma \ref{properties_stoch_sys}, there exist
  $\lambda\ge 0$, $V_1\ge 0$, such that $C_rV_1=0$, i.e.,
  $\langle C_r^\top C_r, V_1\rangle_F=0$ and
  \begin{align*}
     \left(\mathcal{L}_{A_r}+\Pi_{N_{r}}\right)^*(V_1)&=\lambda V_1\;.
  \end{align*}
The scalar products of \eqref{eq:Pgram_trunc},
  \eqref{eq:Qgram_trunc} with $V_1$ yield
  \begin{align}\label{eq:Pgram_truncV}
    \lambda\langle \Sigma_r^{-1},V_1\rangle_F &\le -\langle
                                              \Sigma_r^{-1}B_rB_r^\top \Sigma_r^{-1},V_1\rangle_F-\langle
                                              \Pi_{M_{r}}\left(\Sigma_{2, n-r}^{-1}\right),V_1\rangle_F \le 0,\\\label{eq:Qgram_truncV}
    \lambda \langle \Sigma_r,V_1\rangle_F &= \langle
                                          \Sigma_rB_rB_r^\top \Sigma_r,V_1\rangle_F-\langle
                                          \Pi_{M_{r}}\left(\Sigma_{2, n-r}\right),V_1\rangle_F.
  \end{align}
  From the inequality \eqref{eq:Pgram_truncV} it follows that
  $\lambda\le 0$, i.e., $\lambda=0$.\\
  Without loss of generality, let us assume that
  $\sigma_{r+1}=\max\{\sigma_{r+1},\ldots,\sigma_n\}$. Then,
  $\frac{\sigma_{r+1}^2}{\sigma_j}\ge\sigma_{r+1}\ge \sigma_j$ for
  $j=r+1,\ldots,n$, i.e., $\Upsilon=\sigma_{r+1}^2\Sigma_{2, n-r}^{-1}-
  \Sigma_{2, n-r}\ge 0$.
  Subtracting \eqref{eq:Pgram_truncV} multiplied with $\sigma_{r+1}^2$
  from \eqref{eq:Qgram_truncV} we obtain
  \begin{align*}
    0&
       \ge \langle
       \Sigma_rB_rB_r^\top \Sigma_r+\sigma_{r+1}^2\Sigma_r^{-1}B_rB_r^\top \Sigma_r^{-1},V_1\rangle_F
       +\langle \Pi_{M_{r}}(\Upsilon),V_1 \rangle_F\ge 0\;.
  \end{align*}
  In particular, it holds that $B_rB_r^\top \Sigma_rV_1=0$ and therefore
  \begin{align*}
    0&=A_rV_1+V_1A_r^\top +\Pi_{N_{r}}(V_1)
       =(A_r-B_rB_r^\top \Sigma_r)V_1+V_1(A_r-B_rB_r^\top \Sigma_r)^\top +\Pi_{N_{r}}(V_1) \;,
  \end{align*}
  contradicting the stability of the reduced closed-loop system
 by Theorem \ref{thm:red_stable_bound}.
\end{proof}

\subsection{Reduced order controller}
  Given a reduced model of an unstable system, it is a natural idea
to use it for stabilization. This has been discussed in
\cite{lqg_balancing} for deterministic systems. In the stochastic setup,
the problem is even more involved, and we just sketch some questions.\\
Consider again the systems \eqref{full_system} and \eqref{red_system}  given
by the data \eqref{part_bal}. We partition the balancing transformation matrix
 as $S_b=\left[
  \begin{smallmatrix}
    S_{b,r}^\top\\ \star
  \end{smallmatrix}
\right]$, where  $S_{b,r}^\top$ contains the first $r$ rows. The state
$x_r$ of the
reduced system \eqref{red_system} approximately satisfies $ x_r=S_{b,r}^\top
x$. If a state feedback control $u=F_rx_r$ stabilizes \eqref{red_system},
i.e., $(A_r-B_rF_r,N_{i, r})$ is stable, then we may choose $u=F_r S_{b,r}^\top
x$ as a candidate to stabilize the original system. By Theorem \ref{thm:red_stable_bound} we
can try $F_r=-B_r^\top\Sigma_r$. This choice is also natural as the LQG reduced systems is designed based on negelecting unimportant information in the original stabilizing feedback control, see Proposition \ref{energy_est} (b). For that reason, the truncated singular values $\sigma_{r+1}, \dots, \sigma_n$ are a good indicator for the stabilization by the reduced feedback. Unfortunately, we cannot give detailed a-priori estimates for suitable $r$, such that $u=-B_r^\top\Sigma_r S_{b,r}^\top x$ stabilizes
\eqref{full_system}. But, of course, we can check the closed-loop a
posteriori for stability. This will be done in an example in Section
\ref{sec_num}. \\
Pursuing the idea further, we may also try to design a reduced dynamic
compensator for \eqref{full_system}. In our setup, this could proceed
via the reduced observer system
\begin{align}\label{eq:reduced_observer_system}
  dx_r(t) &= \left[A_r x_r(t) +  {B_r} u(t)+K_r(C_rx_r(t)-y(t))\right]\,dt + \sum_{i=1}^{q} {N}_{i, r} x_r(t) \,dW_i(t).
\end{align}
Setting $K_r=-\Sigma_r C_r^\top$ and $u=-B_r^\top\Sigma_r x_r$, the closed-loop
system can be shown to be stable for $r=n$. For smaller $r$, stability
may be checked a-posteriori. But there is a more serious problem with
this approach, since the noise terms $dW_i$ usually cannot be
reproduced in the observer. Therefore, a thorough analysis would have
to consider only the deterministic part of
\eqref{eq:reduced_observer_system}. We have not carried out any such
work yet which is part of future studies.

\section{Error analysis and its discussion}\label{sec_error}

In this section, we 
begin with an overview on how the error analysis of LQG balancing is conducted in the deterministic case and address difficulties in using the same techniques in the stochastic setting. Subsequently, we provide error bounds for stochastic LQG balancing and show links to the deterministic gap metric analysis.

\paragraph{Deterministic case ($N_i=0$ and deterministic control $u$)}

Given that $N_i=0$, the error analysis between \eqref{full_system} and \eqref{red_system} is often conducted in the frequency domain. To do so, one applies the Laplace transformation to \eqref{full_system} and hence obtains  $\mathbf y = \mathbf G \mathbf u$, where 
$\mathbf u$, $\mathbf y$ are the Laplace transforms of the input and the  output, respectively, and  $\mathbf G$ is the matrix-valued transfer function of the system. The difference between the full and the reduced system can now be measured based on  $\mathbf G - \mathbf G_r$, where $\mathbf G_r$ is the reduced transfer function. A possible error norm can be the $\mathcal H_\infty$-norm defined by
\begin{align}\label{defineHinfty}
\left\| \mathbf G\right\|_{\mathcal H_\infty} := \sup_{w\in\mathbb R} \left\| \mathbf G(\mathrm i w) \right\|_2 = \sup_{u\neq 0} \frac{\left\| y \right\|_{L^2}}{\left\| u \right\|_{L^2}}.                                                                                                              
\end{align}
However, a more suitable error measure in the LQG balancing context is the so-called gap metric. An error analysis for different types of deterministic settings in this metric can be found in \cite{breiten_morandin_schulze, Curtain_LQG, moeckel_reis_stykel}. A possible definition of the gap metric relies on a normalized (right) coprime factorization of the transfer function, i.e., $\mathbf G(s) = \mathbf N(s) \mathbf M(s)^{-1}$. We refer to \cite{breiten_morandin_schulze, DRV, moeckel_reis_stykel} for more details on this factorization. The normalized coprime factors $\mathbf M, \mathbf N$ can now be used to define the gap metric \cite{defn_gap_metric}:\begin{align*}
\delta(\mathbf G, \mathbf G_r):=\max\left\{\inf_{\Pi\in\mathcal H_\infty}\left\| \begin{bmatrix}     \mathbf M_r\\\mathbf N_r   \end{bmatrix}-\begin{bmatrix}     \mathbf M\\\mathbf N   \end{bmatrix} \Pi\right\|_{\mathcal H_\infty},  \inf_{\Pi\in\mathcal H_\infty} \left\| \begin{bmatrix}     \mathbf M\\\mathbf N   \end{bmatrix}-\begin{bmatrix}     \mathbf M_r\\\mathbf N_r   \end{bmatrix} \Pi\right\|_{\mathcal H_\infty}
     \right\}.                                                                                                                                                                                                                                                                                                                                                                                                                                                                                                                                                                                                                                                                                                                                                              \end{align*}
A time-domain interpretation of this distance is, e.g., discussed in \cite{Ball_Sasane, DRV, moeckel_reis_stykel}. Given and $L^2$-input-output pair $u$ and $y$, the gap metric guarantees the existence of a reduced $L^2$-pair $u_r$ and $y_r$, so that we have \begin{align}\label{bound_based_on_gap}
 \left\|\begin{bmatrix} u-u_r\\ y-y_r\end{bmatrix}\right\|_{L^2}
 \leq \delta(\mathbf G, \mathbf G_r)   \left\|\begin{bmatrix} u\\ y\end{bmatrix}\right\|_{L^2}.                                                                                                                                                 \end{align}
A bound for the gap metric is often found using the following estimate 
\begin{align}\label{compute_gap}
\delta(\mathbf G, \mathbf G_r)\leq  \left\| \begin{bmatrix}     \mathbf M\\\mathbf N   \end{bmatrix}-\begin{bmatrix}     \mathbf M_r\\\mathbf N_r   \end{bmatrix} \right\|_{\mathcal H_\infty}.                                                                                                                                                                                                                                                                                                                                                                                                                                                                                                                                                                                                                                                                                                                                                              \end{align}
The $\mathcal H_\infty$-error in \eqref{compute_gap} can be determined based on the time-domain representation of this norm given in \eqref{defineHinfty}. This means, that we can work with system realizations of the transfer functions $\begin{bmatrix}     \mathbf M\\\mathbf N   \end{bmatrix}$, $\begin{bmatrix}     \mathbf M_r\\\mathbf N_r   \end{bmatrix}$ and compute the $L^2$-distance of two associated systems in order to find a bound for the gap metric. However, working with stochastic systems causes various issues since frequency-domain considerations can not be applied. This is due to the fact that the ``derivatives'' in \eqref{full_system} are no longer classical functions not allowing for a Laplace transformation. Therefore, a gap metric study is not possible but our error analysis will rely on generalized system realizations of normalized coprime factorizations. In particular, a reduced input-output pair is supposed to be constructed, so that we find an estimate of the form given in \eqref{bound_based_on_gap}.

\paragraph{Stochastic error analysis}
In order to conduct a gap-metric type error analysis, we construct a pair $u_r, y_r$ that is supposed to well approximate $u, y$. In order to show the error between both vectors, system \eqref{full_system} is rewritten. To be more precise, its
 input-output pair can be parameterized as 
\begin{equation}
 \begin{aligned}\label{sys_par}
 dx(t) =& [\bar A x(t) +  B v(t)]dt + \sum_{i=1}^{q}  N_i x(t) dW_i(t),\\
 \bar y(t):=& \begin{bmatrix} u(t)\\y(t)\end{bmatrix}= \bar C x(t) + \begin{bmatrix} v(t)\\0\end{bmatrix},\quad t\geq 0,
\end{aligned}
\end{equation}
where $\bar A = A- B B^\top Q$, $\bar C= \begin{bmatrix} - B^\top Q \\ C\end{bmatrix}
$ and   $v(t) = B^\top Q x(t) + u(t)$. We can interpret \eqref{sys_par} as a generalized realization (additional $dW_i$ terms) of the coprime factors  $\begin{bmatrix}     \mathbf M\\\mathbf N   \end{bmatrix}$.
In some way, \eqref{sys_par} mimics an asymptotically mean square stable control system since an open-loop system with coefficients $(\bar A, N_i)$ is asymptotically mean square stable due to Theorem \ref{thm:StochLQ}. However, $v$ depends on the solution itself besides depending on $u$. On the other hand, $\bar y$ represents and input-output pair rather than an output. If $N_i = 0$, $v$ and \eqref{sys_par} are called driving-variable and driving-variable system, respectively. The relation between such driving-variable and input-output systems are nicely described in \cite{guiver_opmeer}.\smallskip

We investigate a particular input-output pair of the reduced system fixing control $u_r(t) = B_r^\top \Sigma_r x_r(t)+ B^\top Q x(t)+u(t)$ ($\Sigma_r = \diag(\sigma_1, \dots, \sigma_r)$), since this allows to rewrite the reduced model as \begin{equation}
 \begin{aligned}\label{red_par}
 d x_r(t) =& [\bar A_{r} x_r(t) +  B_r v(t)]dt + \sum_{i=1}^{q} N_{i, r} x_r(t) dW_i(t),\\
 \bar y_{r}(t):=& \begin{bmatrix} u_r(t)\\ y_r(t)\end{bmatrix}= \bar C_{r} x_r(t) + \begin{bmatrix} v(t)\\0\end{bmatrix},\quad t\geq 0,
\end{aligned}
\end{equation}
setting $\bar A_{r} =  A_r- B_r B_r^\top \Sigma_r$ and $\bar C_{r}= \begin{bmatrix} -  B_r^\top \Sigma_r \\  C_r\end{bmatrix}$. Again, \eqref{red_par} can be interpreted as generalized driving variable system or system realization of the reduced coprime factorization. The following theorem establishes an error between the original pair $\begin{bmatrix}  u\\ y \end{bmatrix} $ and the chosen reduced pair $\begin{bmatrix} u_r\\ y_r\end{bmatrix} $. The result relies on $L^2_T$-error estimates between \eqref{sys_par} and \eqref{red_par}.
\begin{thm}\label{main_error_bound}
 Let $u_r(t) = - B_r^\top \Sigma_r x_r(t)+ B^\top Q x(t)+u(t)$ and $y_r$ the reduced order output associated to this input. Given $x_0=0$ and $x_{0, r}=0$, we have \begin{align}\label{bound_finite_time}
 \left\|\begin{bmatrix} u-u_r\\ y-y_r\end{bmatrix}\right\|_{L^2_T}
 \leq 2\sum_{k=r+1}^n \frac{\sigma_k}{\sqrt{1+\sigma_k^2}}    \Bigg(\left\|\begin{bmatrix} u\\ y\end{bmatrix}\right\|_{L^2_T}^2 + \mathbb E\left[x(T)^\top Q x(T)\right] \Bigg)^{\frac{1}{2}}.                                                                                                                                                             \end{align}
If it, furthermore, holds that the input and the state are square integrable on $\Omega\times [0, \infty)$, i.e., $u, x\in L^2$, then we have \begin{align}\label{bound_infinite_time}
 \left\|\begin{bmatrix} u-u_r\\ y-y_r\end{bmatrix}\right\|_{L^2}
 \leq 2\sum_{k=r+1}^n \frac{\sigma_k}{\sqrt{1+\sigma_k^2}}   \left\|\begin{bmatrix} u\\ y\end{bmatrix}\right\|_{L^2}.                                                                                                                                                 \end{align}
\end{thm}
\begin{proof}
We improve the readability of this paper by moving the proof to Appendix \ref{appB}.
\end{proof}     
As a consequence of Theorem \ref{main_error_bound}, we observe that the singular values $\sigma_k$ deliver a good a-priori criterion for the choice of $r$ because removing only small singular values leads to a small bound for the error between the original and the reduced input-output pair. However, this argument is only valid if the (finite time) cost functional and, in case of \eqref{bound_finite_time}, the terminal value $x(T)$ is not too large. The result in \eqref{bound_infinite_time} is a gap-metric type estimate in the sense of \eqref{bound_based_on_gap}. 
We formulate a special case of Theorem \ref{main_error_bound} for $u$ being a stabilizing feedback control. 
\begin{kor}\label{cor_feedback}
 Let $u(t) = - B^\top Q x(t)+u^{(1)}(t)$ and  $u_r(t) = - B_r^\top \Sigma_r x_r(t)+u^{(1)}(t)$ with $u^{(1)}\in L^2_T$. Given $x_0=0$ and $x_{0, r}=0$, we have \begin{align*}
 \left\|\begin{bmatrix} u-u_r\\ y-y_r\end{bmatrix}\right\|_{L^2_T}
 \leq 2\sum_{k=r+1}^n \frac{\sigma_k}{\sqrt{1+\sigma_k^2}}  \left\|u^{(1)}\right\|_{L^2_T}.                                                                                                                                                 \end{align*}
\end{kor}
\begin{proof}
By \eqref{estimate_on_v}, we have $\left\|\begin{bmatrix} u\\ y\end{bmatrix}\right\|_{L^2_T}^2 + \mathbb E\left[x(T)^\top Q x(T)\right] = \left\|B^\top Q x+u\right\|_{L^2_T}^2 = \left\|u^{(1)}\right\|_{L^2_T}^2$. For that reason, this result is a direct consequence of Theorem \ref{main_error_bound}.
\end{proof}   
Corollary \ref{cor_feedback} tells that the stabilizing feedback control $u = - B^\top Q x+u^{(1)}$ and its output can be well-approximated by the reduced feedback $u_r = - B_r^\top \Sigma_r x_r+u^{(1)}$ and the associated output in case the truncated singular values are small.  
\smallskip

We draw our attention back to open-loop controls and discuss the benefit of Theorem \ref{main_error_bound} in this context since this might not be obvious seeing that $u_r$ depends on the original state $x$. Therefore, it seems that we did not gain much from the practical point of view although we found a good candidate for an approximating input-output pair. However, there is a fundamental difference between stochastic and deterministic settings since in the context of stochastic differential equations, there are many problems that cannot be solved in moderate high dimensions $n$ even though one is willing to simulate the original system \eqref{full_system}. To be more precise, one often needs to compute conditional expectations of the form \begin{align*}
g(x):=\mathbb{E}[f(y(t))| x(s)=x],\quad x\in\mathbb R^n, \quad s< t,                                                                                                                                                                                                                                                        \end{align*}
which is the expectation of some quantity of interest $f(y)$ at time $t$ given that the state at time $s$ is $x$. Such objects occur in stochastic optimal stopping problems, e.g., in the context of pricing (Bermudan) options in finance. In order to find an approximation $g(\cdot) \approx\sum_{k=1}^{K} \widehat{\beta}_k\psi_{k}(\cdot)$ of the unknown function $g$, where $\psi_{1}, \ldots, \psi_{K}$ is some suitable (polynomial) basis, we have to solve the least squares problem
\begin{equation}
\label{eq:least-squares-regression}\widehat{\beta} := \argmin_{\beta
\in\mathbb{R}^{K}} \sum_{i=1}^{M} \left|  f(y(t)^{i}) - \sum_{k=1}^{K} \beta_{k}
\psi_{k}(x(s)^{i}) \right|  ^{2},
\end{equation}
where $y(t)^{i}$ and $x(s)^{i}$ i.i.d.~samples of the random variables $y(t)$ and $x(s)$, respectively. Notice that \eqref{eq:least-squares-regression} is the the discretized version by Monte Carlo of the original continuous problem
$\min_{\beta\in\mathbb{R}^{K}} \mathbb{E}
\left|  f(y(t)) - \sum_{k=1}^{K} \beta_{k} \psi_{k}(x(s)) \right|  ^{2}$. Now, solving the regression problem in \eqref{eq:least-squares-regression} requires a huge computational effort already in moderate high dimensions since regression suffers from the curse of dimensionality. This often makes this procedure infeasible for dimensions $n\geq 10$. Therefore, a possible strategy can be to simulate the original system \eqref{full_system} in order to determine the reduced order input $u_r$ defined in Theorem \ref{main_error_bound} that gives a good approximation $y_r$ of $y$. If $r$ is sufficiently small, one can then solve \eqref{eq:least-squares-regression} in the reduced system, in which the impact of the curse of dimensionality is drastically decreased. This leads to a good estimate $g_r$ (defined on $\mathbb R^r$) of the original $g$.\smallskip

We finally investigate the scenario, in which we do not intent to simulate the original system \eqref{full_system} but an open-loop control $u$ is used. Fortunately, Theorem \ref{main_error_bound} also provides a bound for the distance between $u_r$ (defined within this theorem) and the original input $u$. For that reason, we know that $u$ and $u_r$ must be close if the truncated singular values of the system are small. Subsequently, we can use that the (reduced) output is Lipschitz continuous in the control term. This is proved in the following lemma.
\begin{lem}
Given the reduced order model \eqref{red_system} with $x_{0, r}=0$, then there exists a constant $\gamma_T>0$ such that 
\begin{align}\label{y_lipschitz}
 \left\|y_r(\cdot, 0, u_r)\right\|_{L^2_T}\leq \gamma_T  \left\|u_r\right\|_{L^2_T}
\end{align}
for all $u_r\in L^2_T$.
\end{lem}
\begin{proof}
We use equation \eqref{obs_gram} associated to the balanced realization with diagonal solution $\Sigma_n$. We can now exploit the partition in \eqref{part_bal} and  evaluate the left upper block of the balanced version of the matrix equation \eqref{obs_gram}. This yields the following inequality \begin{align}\label{obs_gram_red}
  A_r^\top \Sigma_r + \Sigma_r A_r  + \sum_{i, j=1}^q N_{i, r}^\top \Sigma_r N_{j, r} k_{ij}+ C_r^\top C_r -\Sigma_r B_r B_r^\top \Sigma_r \leq 0,                                                                                                                                                                                                                                                                                                                                                                                                                                                                                                                                                                  \end{align}
where $\Sigma_r= \diag(\sigma_{1},\ldots,\sigma_r)$ contains the first $r$ singular values of the system. Applying \eqref{estimate_quadr_functional} to the reduced system with initial state zero and setting $X=\Sigma_r$, we obtain \begin{align}
\nonumber\mathbb E \left[x_r(t)^\top \Sigma_r x_r(t)\right] 
=&\int_0^t \mathbb E\left[x_r(s)^\top\left(A_r^\top \Sigma_r + \Sigma_r A_r  +\sum_{i, j=1}^q N_{i, r}^\top \Sigma_r N_{j, r} k_{ij} \right) x_r(s)\right] ds\\ \label{estimate_red_sys}&+ 2 \int_0^t \mathbb E\left\langle B_r^\top \Sigma_r x_r(s),  u_r(s) \right\rangle_2ds.
\end{align}
With $2 \left\langle B_r^\top \Sigma_r x_r(s),  u_r(s) \right\rangle_2\leq \left\| B_r^\top \Sigma_r x_r(s) \right\|_2^2 + \left\| u_r(s) \right\|_2^2$ and \eqref{obs_gram_red}, identity \eqref{estimate_red_sys} becomes \begin{align*}
\mathbb E \left[x_r(t)^\top \Sigma_r x_r(t)\right] 
\leq& \left\| u_r \right\|_{L^2_t}^2-\left\| y_r \right\|_{L^2_t}^2  + 2 \int_0^t \mathbb E \left\| B_r^\top \Sigma_r x_r(s) \right\|_2^2 ds \\ \leq & \left\| u_r \right\|_{L^2_t}^2-\left\| y_r \right\|_{L^2_t}^2  + 2 b_r \int_0^t \mathbb E \left[x_r(s)^\top \Sigma_r x_r(s)\right] ds,                                                                                                                                                                                                                                                                                                                                                                                                                                          \end{align*}
where
$b_r:= \left\| B_r^\top \Sigma_r^{\frac{1}{2}} \right\|_2^2$. Gronwall's Lemma \ref{gronwall} for  $t=T$ leads to \begin{align*}
\mathbb E \left[x_r(T)^\top \Sigma_r x_r(T)\right] &\leq\left\| u_r \right\|_{L^2_T}^2-\left\| y_r \right\|_{L^2_T}^2 + \int_0^T (\left\| u_r \right\|_{L^2_s}^2-\left\| y_r \right\|_{L^2_s}^2) 2 b_r \expn^{2 b_r(T-s)}ds\\
&= \int_0^T (\left\| u_r(s) \right\|_{2}^2-\left\| y_r(s) \right\|_{2}^2) \expn^{2 b_r(T-s)}ds 
\end{align*}
using integration by parts in the last step. Therefore, we have \begin{align*}
 \left\|  y_r \right\|_{L^2_T}^2  \leq \int_0^T \left\| y_r(s) \right\|_{2}^2  \expn^{2 b_r(T-s)}ds \leq \int_0^T \left\| u_r(s) \right\|_{2}^2\expn^{2 b_r(T-s)}ds\leq \expn^{2 b_r T} \left\|  u_r \right\|_{L^2_T}^2.                                                                                                                         \end{align*}
This concludes the proof.
\end{proof}
By the linearity of $y_r$ in $u_r$,  \eqref{y_lipschitz} means that controls being close to each other lead to similar outputs. Therefore, only a slight deviation between $y_r(\cdot, 0, u_r)$ and $y_r(\cdot, 0, u)$ is expected. The smallest constant in \eqref{y_lipschitz} is $\gamma_T= \sup_{u_r\in L^2_T\setminus \{0\}} \frac{\left\|y_r(\cdot, 0, u_r)\right\|_{L^2_T}}{\left\|u_r\right\|_{L^2_T}}$. If there is a Lipschitz constant independent of $T$, we can consider $\gamma_T=\gamma= \sup_{u_r\in L^2\setminus \{0\}} \frac{\left\|y_r(\cdot, 0, u_r)\right\|_{L^2}}{\left\|u_r\right\|_{L^2}}$ in \eqref{y_lipschitz} which is the norm of the input-output operator on the entire positive real line. This holds, e.g., if \eqref{red_system} is asymptotically stable \cite{Hinrichsen_Pritchard}. 
We can now formulate the result when $u_r = u$ is chosen in Theorem \ref{main_error_bound}.
\begin{kor}\label{cor_open_loop}
 Let $x_0=0$, $x_{0, r}=0$, $u\in L^2_T$ and $y_r=y_r(\cdot, 0, u)$. If $\gamma_T>0$ is a constant satisfying \eqref{y_lipschitz}, we have \begin{align*}
 \left\| y-y_r\right\|_{L^2_T}
 \leq 2 (1+\gamma_T)\sum_{k=r+1}^n \frac{\sigma_k}{\sqrt{1+\sigma_k^2}}    \Bigg(\left\|\begin{bmatrix} u\\ y\end{bmatrix}\right\|_{L^2_T}^2 + \mathbb E\left[x(T)^\top Q x(T)\right] \Bigg)^{\frac{1}{2}}.                                                                                                                                                             \end{align*}
If additionally holds that $u, x\in L^2$, then we have \begin{align*}
 \left\|y-y_r\right\|_{L^2_T}
 \leq 2(1+\gamma_T)\sum_{k=r+1}^n \frac{\sigma_k}{\sqrt{1+\sigma_k^2}}   \left\|\begin{bmatrix} u\\ y\end{bmatrix}\right\|_{L^2}.                                                                                                                                                 \end{align*}
\end{kor}
\begin{proof}
 It holds that \begin{align}\nonumber
       \left\| y(\cdot, 0, u)-y_r(\cdot, 0, u)\right\|_{L^2_T} &\leq 
       \left\| y(\cdot, 0, u)-y_r(\cdot, 0, u_r)\right\|_{L^2_T}+\left\| y_r(\cdot, 0, u_r)-y_r(\cdot, 0, u)\right\|_{L^2_T}\\ \label{triangle_est} &\leq \left\| y(\cdot, 0, u)-y_r(\cdot, 0, u_r)\right\|_{L^2_T}+\gamma_T\left\|  u_r-u\right\|_{L^2_T},
               \end{align}
where $u_r$ is defined as in Theorem \ref{main_error_bound}. Applying \eqref{bound_finite_time} to both terms in \eqref{triangle_est} yields the first estimate. If $u, x\in L^2$ holds, we can use \eqref{bound_infinite_time} instead and obtain the second inequality. This concludes the proof.
\end{proof}
According to Theorem \ref{main_error_bound} the singular values of \eqref{full_system} can be used a-priori to find a suitable dimension $r$ of an accurate reduced system \eqref{red_system} using a control that is possibly not available. Corollary \ref{cor_open_loop} now additionally tells us that this unavailable control can be replaced by the original system control if the norm of the input-output operator is not too large. Computing this norm $\gamma_T$ is feasible in small dimensions $r$ without causing large computation cost. Hence, $\gamma_T$ provides an a-posteriori criterion for a good approximation of $y(\cdot, 0, u)$ by $y_r(\cdot, 0, u)$. We finally provide a bound that neither contains the state $x$ nor the output $y$ of the original system.
\begin{thm}\label{main_error_bound2}
 Let $y=y(\cdot, 0, u)$ and $y_r = y_r(\cdot, 0, u)$ and $u\in L^2_T$. Then, 
 \begin{align*}
\left( \mathbb E\int_0^T \expn^{-\beta t}\left\|y(t)-y_r(t)\right\|_{2}^2 dt \right)^{\frac{1}{2}}
 \leq 2\sum_{k=r+1}^n \sigma_k   \left( \mathbb E\int_0^T \expn^{-\beta t}\left\|u(t)\right\|_{2}^2 dt \right)^{\frac{1}{2}},                                                                                                                             \end{align*}
 where $\beta =\max\{\left\|B^\top Q^{\frac{1}{2}}\right\|_2^2,  \left\|C P^{\frac{1}{2}}\right\|_2^2\}$.
\end{thm}
\begin{proof}
We present the proof in Appendix \ref{appC}.
\end{proof}
The bound of Theorem \ref{main_error_bound2} is practically computable since it does not involve variables of the original system \eqref{full_system}. However, a high accuracy cannot be expected since it is a worst-case bound (based on Gronwall's lemma) that also captures systems with exponentially growing states that might not be approximated well with the underlying dimension reduction scheme.  Therefore, the result of Theorem \ref{main_error_bound2} can also be read as a warning that LQG balancing is not working well for all types of unstable open-loop systems (satisfying our assumptions) even though the truncated singular values are small.
                                                            
\section{Numerical examples}\label{sec_num}

For $t\in [0, T]$, we consider the following $2$D stochastic heat equation with
Neumann boundary conditions and scalar noise ($q=1$):
\begin{align*}
\frac{\partial X(t, \zeta)}{\partial t} &= \alpha\Delta X(t, \zeta)+f(\zeta) u(t)+
 \nu g(\zeta)X(t, \zeta)\frac{\partial W(t)}{\partial t},\quad \zeta\in [0, \pi]^2,\\
\frac{\partial X(t, \zeta)}{\partial\mathbf{n}} &= 0,\quad \zeta\in \partial [0, \pi]^2,\quad X(0, \zeta)\equiv 0,
\end{align*}
where $\alpha,\nu>0$ and $f, g$ are bounded functions on $[0, \pi]^2$. 
We set $\alpha=0.2$, $\nu=2$, $f(\zeta) = 1_{[\frac{\pi}{4}, \frac{3 \pi}{4}]^2}(\zeta)$, $g(\zeta) = \expn^{-\left\vert \zeta_1-\frac{\pi}{2}\right\vert-\zeta_2}$ and 
$H=L^2([0, \pi]^2)$ to be the solution space for the mild solution of the stochastic partial differential equation (SPDE). In this contex, let $\langle\cdot, \cdot \rangle_H$ denote the inner product in $H$ and $\left\|\cdot\right\|_H$ the corresponding norm. The output is the mean temperature on the
uncontrolled area,
\begin{align*}
Y(t)= \mathcal C X(t, \zeta):= \frac{4}{3 \pi^2}\int_{[0, \pi]^2\setminus [\frac{\pi}{4}, \frac{3 \pi}{4}]^2} X(t, \zeta) d\zeta.
\end{align*}
We discretize this SPDE by a spectral Galerkin method according to \cite{redmannbenner}. 
The eigenvalues of the Neumann Laplacian on $[0, \pi]^2$ are given by
$\lambda_{i j}=-(i^2+j^2)$ and the corresponding eigenvectors
representing an ONB of $H$ are $h_{i
  j}=\frac{f_{ij}}{\left\|f_{i j}\right\|_H}$, where $f_{i
  j}=\cos(i\cdot) \cos(j\cdot)$. We order these eigenvalues and write
$\lambda_k$ and $h_k$ for the $k$-th largest eigenvalue and the
associated eigenvector, respectively. We obtain a system of the form
\eqref{full_system} with matrices $C^\top=\left(\mathcal C
  h_k\right)_{k=1, \ldots, n}$, $A=\alpha\diag(\lambda_1, \lambda_2,
\dots)=\alpha\diag(0, -1, \dots) $,  $N_1 =\nu\left(\left\langle g h_i,
    h_k\right\rangle_H\right)_{k, i=1, \ldots, n}$,
$B=\left(\left\langle f , h_k\right\rangle_H\right)_{k=1, \ldots, n}$. We observe that this spatial discretization is unstable but the requirements for applying LQG balancing are fulfilled.

Gramians $Q$ and $P$ according to their definitions in
Theorem \ref{thm:StochLQ} and in \eqref{reach_gram} can now be computed.
For $Q$ we have used a fixed point iteration with $Q_0=I$ and
$Q_{k+1}$ being the stabilizing solution of the Riccati equation
\begin{align*}
  \mathcal L_{A}(Q_{k+1})+\Pi_N(Q_k)+C^\top C-Q_{k+1}BB^\top Q_{k+1}=0.
\end{align*}
This converges quite fast to the Gramian $Q$,
e.g., \cite[Sec.~5.4.3]{damm}. The Gramian $P$ is computationally more
involved. By our error analysis in Section \ref{sec_error}, it is natural to seek for a $P$ with a large number of small eigenvalues, so that we aim to find the Gramian with minimal trace subject to \eqref{reach_gram}. However, we do not have a linear matrix inequality (LMI) formulation for $P$ but rather for its inverse. Therefore, we have rewritten \eqref{reach_gram} as the LMI
\begin{align*}
  \left[
  \begin{array}{cc}
A^\top P^{-1}+P^{-1}A+N_1^\top P^{-1}N_1-C^\top C&P^{-1}B\\ B^\top P^{-1}&-I
  \end{array}
\right]\le 0\;,\quad P^{-1}\ge 0
\end{align*}
and maximized the trace of $P^{-1}$ using the solver Mosek \cite{mosek}
with the Matlab package Yalmip \cite{Lofberg2004}. However, this might not ensure the same approximation quality as when being able to solve for $P$ directly. For the computation of $P^{-1}$ in dimension $n=100$, it
took about $5$ to $6$ minutes, where the empirical complexity is about
$n^6$. Therefore, we did not consider larger systems.

\begin{figure}[h]\centering
    \begin{minipage}{.48\linewidth}
%
%
\definecolor{mycolor1}{rgb}{0.00000,0.44700,0.74100}%
\begin{tikzpicture}

\begin{axis}[%
width=5.2cm,
height=3.7cm,
at={(0cm,0cm)},
scale only axis,
xmin=0,
xmax=100,
ymode=log,
ymin=1e-10,
ymax=100000000,
yminorticks=true,
axis background/.style={fill=white},
legend style={legend cell align=left, align=left, draw=white!15!black}
]
\addplot [color=blue, draw=none, mark=asterisk, mark options={solid, blue}]
  table[row sep=crcr]{%
1	1033365.29624788\\
2	25930.7287015434\\
3	6344.13620845148\\
4	3339.79769710032\\
5	257.37668799938\\
6	162.075390624215\\
7	44.8023854880781\\
8	29.5097385748028\\
9	22.8172000131204\\
10	10.5164298708629\\
11	6.17739881428782\\
12	2.15218450536298\\
13	1.17851454108732\\
14	0.650614087683786\\
15	0.386366883552092\\
16	0.258748542694428\\
17	0.123513549801187\\
18	0.0715605346662529\\
19	0.0427968020017345\\
20	0.0226446779488288\\
21	0.0200917260301841\\
22	0.0127930525197452\\
23	0.0105039784315886\\
24	0.00856336890708062\\
25	0.0031593264185311\\
26	0.00183730071327713\\
27	0.00149600971074982\\
28	0.00129531583591949\\
29	0.00081829588373124\\
30	0.00056574293729731\\
31	0.000269119647446794\\
32	0.000193957329606056\\
33	0.000164873020418211\\
34	0.000136903209449554\\
35	8.926693279874e-05\\
36	7.71641664861203e-05\\
37	3.93658471613572e-05\\
38	3.22730425875754e-05\\
39	1.39233920826911e-05\\
40	1.10467140539788e-05\\
41	7.70804993513795e-06\\
42	4.72584598898066e-06\\
43	4.67798116664798e-06\\
44	3.50664070054633e-06\\
45	2.96602701840333e-06\\
46	1.71046220524656e-06\\
47	1.44495187469917e-06\\
48	1.26376614236687e-06\\
49	1.21928979620241e-06\\
50	6.52041451909757e-07\\
51	4.77236921627765e-07\\
52	3.06520657271979e-07\\
53	2.83146377510135e-07\\
54	7.99894115100328e-08\\
55	7.9216377987755e-08\\
56	7.86210091567509e-08\\
57	7.86197870956311e-08\\
58	7.75942426081739e-08\\
59	7.69300161639186e-08\\
60	7.67038359716717e-08\\
61	7.50020949933507e-08\\
62	7.49988912938304e-08\\
63	7.47352054646945e-08\\
64	7.36131986266173e-08\\
65	7.36046858249319e-08\\
66	7.27030927922659e-08\\
67	7.26782606035337e-08\\
68	7.23625942362352e-08\\
69	6.98440392054749e-08\\
70	6.98218005801136e-08\\
71	6.94306197438703e-08\\
72	6.5783017301931e-08\\
73	6.57575123438141e-08\\
74	6.34505633189952e-08\\
75	6.14664201264342e-08\\
76	6.14549589800439e-08\\
77	5.79985487211252e-08\\
78	5.58950055452651e-08\\
79	5.58675075691122e-08\\
80	5.58266673853582e-08\\
81	5.51144520338425e-08\\
82	5.21358400549338e-08\\
83	4.90297025948307e-08\\
84	4.59064829258633e-08\\
85	4.34776765442793e-08\\
86	4.34327288103325e-08\\
87	3.93407146775226e-08\\
88	3.69454413971638e-08\\
89	3.6888235100547e-08\\
90	3.61311529960724e-08\\
91	3.610165096754e-08\\
92	3.10639845590325e-08\\
93	3.04621361987651e-08\\
94	2.82074278403812e-08\\
95	2.63187796549096e-08\\
96	2.60356204517635e-08\\
97	2.45822209810598e-08\\
98	1.73776777115788e-08\\
99	9.38815844790724e-09\\
100	1.89366096688008e-09\\
};
\addlegendentry{Singular values $\sigma_j$}

\addplot [color=mycolor1]
  table[row sep=crcr]{%
11	1.89366096688008e-09\\
11	1033365.29624788\\
};
\addlegendentry{$\sigma{}_{\text{r+1}}\text{/}\sigma{}_\text{1}\text{=1.0177e-05}$}

\addplot [color=mycolor1]
  table[row sep=crcr]{%
0	10.5164298708629\\
100	10.5164298708629\\
};

\end{axis}

\end{tikzpicture}%
    \end{minipage}
\caption{Decay singular values of discretized heat equation for $n=100$.}\label{fig:HSV}
\end{figure}
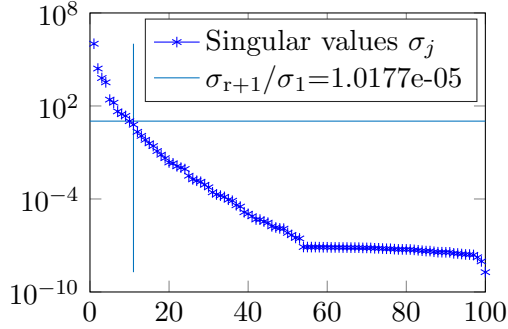

\begin{figure}[h]
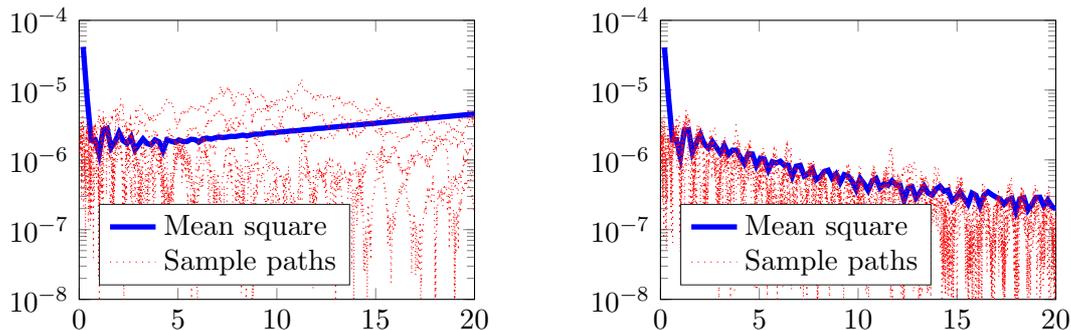
\centering
    \begin{minipage}{.48\linewidth}
      \input{./error_open-loop20a.tex}
    \end{minipage} \hfill
  \begin{minipage}{.48\linewidth}
       \input{./error_closed-loop20a.tex}
      \end{minipage}
\caption{Open-loop and closed-loop error systems driven by $L^2$-input $u$.}\label{fig:errCL}
\end{figure}

For $n=100$, we choose the reduced order $r=10$. The decay of the
singular values $\sigma_j$ is shown in Fig. \ref{fig:HSV}. As pointed out in Section \ref{sec_error}, the truncated singular values determine the error of  LQG balanced truncation. We observe that $r=10$ provides very small $\sigma_{r+1}, \dots, \sigma_n$ relative to $\sigma_1$.
We combine \eqref{full_system} and \eqref{red_system} ($u_r=u$) with zero initial states and define $\xi = \left[
        \begin{smallmatrix}
          x\\x_r
        \end{smallmatrix}
\right]$ leading the open-loop error system  
 \begin{align*}
  d\xi&= \left(\left[
        \begin{smallmatrix}
          A&0\\0&A_r
        \end{smallmatrix}
\right]\xi+ \left[
        \begin{smallmatrix}
          B\\B_r
        \end{smallmatrix}
\right]u\right)\,dt+\left[
        \begin{smallmatrix}
          N_1&0\\0&N_{1, r}
        \end{smallmatrix}
\right]\xi\,dW,\quad y_\xi=\left[
        \begin{smallmatrix}
          C&-C_r
        \end{smallmatrix}
\right]\xi,
 \end{align*}
where $y_\xi= y-y_r$. Secondly, we introduce a closed-loop version by  
 \begin{align*}
  d\xi&= \left(\left[
        \begin{smallmatrix}
          A-BB^\top Q &0\\0&A_r-B_rB_r^\top\Sigma_r
        \end{smallmatrix}
\right]\xi+ \left[
        \begin{smallmatrix}
          B\\B_r
        \end{smallmatrix}
\right]u\right)\,dt+\left[
        \begin{smallmatrix}
          N_1&0\\0&N_{1, r}
        \end{smallmatrix}
\right]\xi\,dW,\quad y_\xi=\left[
        \begin{smallmatrix}
          C&-C_r
        \end{smallmatrix}
\right]\xi,
 \end{align*}
i.e., a stabilizing feedback control is used.  We have computed $t\mapsto\sqrt{\mathbb E \left \vert y_\xi(t)\right\vert^2}$ (blue graphs) and 
five trajectories $t\mapsto\left \vert y_\xi(t, \omega)\right\vert$ (red dotted graphs) in Fig.\ \ref{fig:errCL}. Notice that the open-loop case is depicted left and the closed-loop scenario is given in the right picture. In both cases, we have used zero initial states and the
 $L^2$-input $u=\frac{\cos(5t)}{t+1}$. 
 
 The mean square error in blue has been computed from a deterministic Lyapunov type ordinary differential equation and the five sample output paths from a drift implicit Euler-Maruyama method.

\begin{figure}[h]
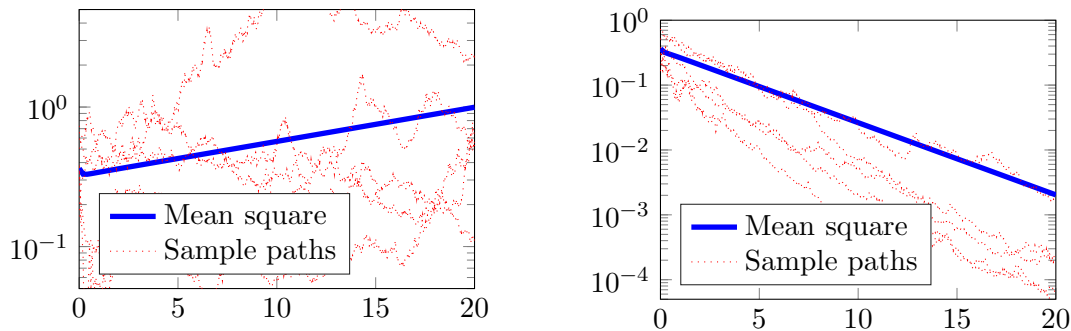
\centering
    \begin{minipage}{.48\linewidth}
      \input{./full_open-loop20a.tex}
    \end{minipage} \hfill
  \begin{minipage}{.48\linewidth}
       \input{./red_contr20a.tex}
     \end{minipage}
\caption{No control vs.\ reduced controller in \eqref{full_system} with random initial state.}
\label{fig:uncontrolled}
\end{figure}

We can see that the error is small in both cases, but  by stability it
decays only in the closed-loop case.
Furthermore, we observed in this example that the reduced feedback controller also stabilizes the
full system. To visualize this effect, we have computed $t\mapsto\sqrt{\mathbb E \left \vert y(t)\right\vert^2}$ (blue graphs) for system \eqref{full_system} with $u\equiv 0$ (Fig.\ \ref{fig:uncontrolled} left), $u(t)=B_r^\top\Sigma_r S_{b,r}^\top x(t)$ (Fig.\ \ref{fig:uncontrolled} right) and a randomly generated initial state $x_0$, where  $S_{b,r}^\top$ are the first $r$ rows of the balancing transformation $S_b$ in Section \ref{sec:balanced-truncation}. As mean square stability is stronger than path-wise stability in the linear case, we see the same effect for the trajectories $t\mapsto\left \vert y(t, \omega)\right\vert$ (red dotted graphs) in Fig.\ \ref{fig:uncontrolled}.

\section{Conclusions}

In this paper, we considered dimension reduction techniques for large-scale stochastic systems. Such schemes are vital in both control and probabilistic settings as many system evaluations are required. In this context, one can think of aiming to investigate statistical properties by sampling methods or the optimal control of spatially discretized stochastic partial differential equations.  These fit into our framework, in which we have studied potentially unstable stochastic differential equations. They occur, for instance if the driving noise is large. Therefore, we have made an essential contribution since many existing model reduction schemes require certain stability conditions. We have introduced a pair of Gramians that are designed in order to characterize dominant subspaces of the underlying stochastic system. In particular, unnecessary direction in closed-loop dynamics have been identified by these Gramians but we have also pointed out their relevance for open-loop controls. These considerations led to a reduced order model that captures many important features of the original one. We proved that, e.g., stabilizability and detectability are preserved. Our dimension reduction procedure further allowed for a detailed error analysis. Based on the error estimates provided in this work, algebraic a-priori criteria for the approximation quality have been found. These error bounds therefore give a clear guidance on how to fix the reduced order dimension. The effectiveness of our method has been demonstrated by applying it to an unstable stochastic heat equation. This means that an infinite dimensional state dynamics could be approximated by a low-order stochastic system. 

\appendix
\section{Supporting lemmas}\label{appA}
\begin{lem}[Gronwall lemma]\label{gronwall}
Given $T>0$ let $z, \alpha: [0, T]\rightarrow \mathbb R$ and $\beta: [0, T]\rightarrow [0, \infty)$ be continuous functions. 
If \begin{align*}
    z(t)\leq \alpha(t)+\int_0^t \beta(s) z(s) ds,\quad t\in[0, T],
   \end{align*}
then for all $t\in[0, T]$, it holds that \begin{align*}
    z(t)\leq \alpha(t)+\int_0^t \alpha(s)\beta(s) \exp\left(\int_s^t \beta(w)dw\right) ds.
   \end{align*}
\end{lem}
   \begin{proof}
  The result can be shown following the steps in \cite[Proposition 2.1]{gronwalllemma}.
 \end{proof}
\begin{lem}\label{lemstochdiff}
Let $a, b_1, \ldots, b_q$ be $\mathbb R^d$-valued processes, where $a$ is $\left(\mathcal F_t\right)_{t\in [0, T]}$-adapted and almost surely Lebesgue integrable and the functions $b_i$ 
are integrable with respect to the mean zero Wiener process $W=(W_1, \ldots, W_q)^\top$ with covariance matrix $K=\left(k_{ij}\right)_{i, j=1, \ldots, q}$. If the process $x$ is given by \begin{align*}
 dx(t)=a(t) dt+ \sum_{i=1}^v b_i(t)dW_i,\quad t\in [0, T].                                                                  \end{align*}
Then, we have \begin{align*}
 \frac{d}{dt}\mathbb E\left[x(t)^\top x(t)\right]=2 \mathbb E\left[x(t)^\top a(t)\right] + \sum_{i, j=1}^v \mathbb E\left[b_i(t)^\top b_j(t)\right]k_{ij}.                                                                                                                          \end{align*}
\end{lem}                                                                                                                            \begin{proof}
We refer to \cite[Lemma 5.2]{redmannspa2} for a proof of this lemma. 
\end{proof}
Let $x$ now be the solution to \eqref{stochstate}. 
As a direct consequence, we obtain the following identity:
\begin{align}\nonumber
\mathbb E \left[x(t)^\top X x(t)\right] &= x_0^\top X x_0 + 2 \int_0^t \mathbb E\left[x(s)^\top X\left(A x(s) + B u(s) \right)\right]ds\\ \nonumber
&\quad + \int_0^t \sum_{i, j=1}^q \mathbb E\left[x(s)^\top N_{i}^\top X N_{j} x(s) \right]k_{ij} ds\\ \nonumber
&=x_0^\top X x_0+\int_0^t \mathbb E\left[x(s)^\top\left(A^\top X + XA  +\sum_{i, j=1}^q N_{i}^\top X N_{j} k_{ij} \right) x(s)\right] ds\\ \label{estimate_quadr_functional}
&\quad + 2 \int_0^t \mathbb E\left\langle B^\top X x(s),  u(s) \right\rangle_2ds,
\end{align}
where $X\geq 0$ is a semidefinite matrix.

\begin{lem}\label{lemdgl}
Let $W$ be as in Lemma \ref{lemstochdiff} and $A, N_i\in\mathbb R^{k\times k}$ be generic matrices.
Suppose that $b$ is an $\mathbb R^k$-valued and $c_0, \dots, c_q$ are scalar $\left(\mathcal F_t\right)_{t\in [0, T]}$-adapted processes in $L^2_T$. If $x$ is given by
 \begin{align}
 dx(t) = [A x(t) +  b(t)\pm \smat 0 \\ c_0(t)\srix]dt + \sum_{i=1}^{q}[ N_{i} x(t)\pm\smat 0 \\ c_i(t)\srix]  dW_i(t), \quad x(0)=0.
\end{align}
Then, for $t\in[0, T]$, we have \begin{align}\nonumber
\mathbb E\left[x(t)^\top D x(t)\right] 
&=\mathbb E\int_0^t \hspace{-0.25cm} x(s)^\top\hspace{-0.1cm}\left(A^\top D+D A^\top+\sum_{i, j=1}^q N_{i}^\top D N_{j} k_{ij}\right)\hspace{-0.1cm}x(s) + 2 x(s)^\top D b(s) ds\\ \label{productruleapplied2}
& \quad\pm d_k\mathbb E\int_0^t \hspace{-0.15cm} 2 x_{2}(s) c_0(s) + \sum_{i, j=1}^q \left(2 n_i  x(s) \pm c_i(s)\right)c_j(s) k_{ij}  ds,
\end{align}
where $D=\diag(d_1, \dots, d_k)\geq 0$, $n_i$ is the last row of $N_i$ and $x_2$ the last entry of $x$.
\end{lem}
\begin{proof}
Applying Lemma \ref{lemstochdiff}, we find
\begin{align*}
\mathbb E\left[x(t)^\top D x(t)\right] &=2 \int_0^t\mathbb E\left[x(s)^\top D \left(A x(s) + b(s) \pm \smat 0 \\ c_0(s)\srix \right)\right]ds\\ 
&\quad+ \int_0^t\sum_{i, j=1}^q \mathbb E\left[\left(N_{i} x(s) \pm \smat 0 \\ c_i(s)\srix\right)^\top D\left(N_{j} x(s) \pm \smat 0 \\ c_j(s)\srix\right)\right]k_{ij} ds\\ 
&=\mathbb E\int_0^t \hspace{-0.25cm} x(s)^\top \hspace{-0.1cm}\left(A^\top D+D A^\top+\sum_{i, j=1}^q N_{i}^\top D N_{j} k_{ij}\right)\hspace{-0.1cm}x(s) + 2 x(s)^\top D b(s) ds\\
& \quad\pm\mathbb E\int_0^t \hspace{-0.15cm}  2 x(s)^\top D \smat 0 \\ c_0(s)\srix + \sum_{i, j=1}^q \left(2 N_{i} x(s) \pm \smat 0 \\ c_i(s)\srix\right)^\top D \smat 0 \\ c_j(s)\srix k_{ij}  ds.
\end{align*}
We observe that $x(s)^\top D\smat 0 \\ c_0(s)\srix=d_k  x_{2}(s) c_0(s)$ and \begin{align*}
 \left(2 N_{i} x(s) \pm \smat 0 \\ c_i(s)\srix\right)^\top D\smat 0 \\ c_j(s)\srix 
 &= d_k \left(2 n_i  x(s) \pm c_i(s)\right)c_j(s),                                                                            \end{align*}
so that the result follows.
\end{proof}

\section{Proof of Theorem \ref{main_error_bound}}\label{appB}
\begin{proof}[Proof of Theorem \ref{main_error_bound}]
Let $(A_n, B_n, C_n, N_{i, n})$ be the balanced realization of \eqref{full_system} with state variable $x_n$. Let us further introduce $A_{k}$ and $N_{i, k}$ as the left upper $k\times k$ blocks of $A_n$ and $N_{i, n}$. Moreover, suppose that $B_{k}$ and $C_{k}$ are the first $k$ rows of $B_n$ and first $k$ columns of $C_n$, $k= r, \dots, n-1$.  We define
\begin{equation}
 \begin{aligned}\label{bal_par}
 dx_k(t) &= [\bar A_{k} x_k(t) +  B_k v(t)]dt + \sum_{i=1}^{q} N_{i, k} x_k(t) dW_i(t),\\
 \bar y_{k}(t)&= \bar C_{k} x_k(t) + \begin{bmatrix} v(t)\\0\end{bmatrix},\quad t\geq 0,
\end{aligned}
\end{equation}
where $\bar A_{k} = A_k- B_k B_k^\top \Sigma_k$, $\bar C_{k}= \begin{bmatrix} - B_k^\top \Sigma_k \\ C_k\end{bmatrix}
$ and $k=r, \dots, n$. Clearly, $k=r$ yields \eqref{red_par}. On the other hand, $k=n$ provides the input-output parameterization of the balanced version of \eqref{full_system} which can be seen by exploiting $B^\top Q x(t) = B_n^\top \Sigma_n x_n(t)$. Consequently, $\bar y_{n}$ coincides with $\bar y$ in \eqref{sys_par}. Therefore, we have 
\begin{align}\label{remove_HSVs_individually}
 \left\|\begin{bmatrix} u-u_r\\ y-y_r\end{bmatrix}\right\|_{L^2_T} = \left\|\bar y-\bar y_{r}\right\|_{L^2_T} 
 \leq \sum_{i=r+1}^n \left\|\bar y_{k}-\bar y_{k-1}\right\|_{L^2_T}, 
 \end{align}
for which we investigate every summand $ \left\|\bar y_{k}-\bar y_{k-1}\right\|_{L^2_T}$ in the following. 
Exploiting the definitions of  $\bar A_{n}$ and $\bar C_{n}$ the balanced version of \eqref{obs_gram} becomes
\begin{align}\label{aux_obs_gram}
\bar A_{k}^\top \Sigma_k + \Sigma_k \bar A_{k} + \sum_{i, j=1}^q  N_{i, k}^\top \Sigma_k N_{i, k} k_{ij}+ \bar C_{k}^\top \bar C_{k} \leq 0
\end{align}
for $k=n$ (here even the equality holds in \eqref{aux_obs_gram}). The inequalities in \eqref{aux_obs_gram} for $k= r, \dots, n-1$ follow by evaluating the $k\times k$ left upper block of inequality with $k=n$. We now define $L_k= \Sigma_k + \Sigma_k^{-1}$. 
Then, adding the balanced versions of both (in)equalities \eqref{obs_gram} and \eqref{reach_gram} yields
\begin{align}\label{aux_reach_gram}
\bar A_{k}^\top L_k + L_k \bar A_{k} 
+ \sum_{i, j=1}^q N_{i, k}^\top L_k  N_{i, k} k_{ij} 
+L_k B_k B_k^\top L_k\leq 0
\end{align}
for $k=n$. The evaluation of the left upper blocks then provides the results for $k= r, \dots, n-1$. We partition \begin{align}\label{partition2}
x_k= \begin{bmatrix} x_{k,1}\\ 
x_{k, 2} \end{bmatrix}, \quad \bar A_{k}= \begin{bmatrix} {\bar A}_{k-1}&\star\\ 
{a}_{21}&\star\end{bmatrix} ,\quad B_k = \begin{bmatrix} {B}_{k-1}\\ {b}_2\end{bmatrix} ,\quad N_{i, k} =  \begin{bmatrix} {N}_{i, k-1}&\star\\ 
{n}_{i, 21}&{n}_{i, 22}\end{bmatrix} .\end{align}
The variable $x_{k, 2}$ is scalar and we omit the index $k$ in $a_{21}, {n}_{i, 21}\in\mathbb R^{1\times (k-1)}$, $b_2\in \mathbb R^{1\times m}$, ${n}_{i, 22}\in\mathbb R$ in order to simplify the notation. We set \begin{align}\label{xplusminus}
x_-=\begin{bmatrix}  x_{k,1}-x_{k-1}\\ x_{k, 2}\end{bmatrix}, \quad x_+=\begin{bmatrix}  x_{k,1}+x_{k-1} \\ x_{k, 2}\end{bmatrix}           
\end{align}
and obtain from \eqref{bal_par} that \begin{subequations}\label{plus_minus_sys}
 \begin{align}\label{minus}
 dx_-(t) &= [\bar A_{k} x_-(t) +  \smat 0 \\ c_0(t)\srix]dt + \sum_{i=1}^{q}[ N_{i, k} x_-(t)+ \smat 0 \\ c_i(t)\srix] dW_i(t),\\ \label{plus}
 dx_+(t) &= [\bar A_{k} x_+(t) +  2 B_k v(t)-\smat 0 \\ c_0(t)\srix]dt + \sum_{i=1}^{q}[ N_{i, k} x_+(t)-\smat 0 \\ c_i(t)\srix]  dW_i(t),
\end{align}
\end{subequations}
where $c_0(t):=a_{21}x_{k-1}(t)+b_2 v(t)$ and $c_i(t):= n_{i, 21} x_{k-1}(t)$. We apply Lemma \ref{lemdgl} to  \eqref{minus} with $b(t) = 0$ and $D= \Sigma_k$. Moreover, we immediately insert \eqref{aux_obs_gram} into the result of this lemma resulting in 
\begin{align}\label{estimate_used_again}
  \mathbb E\left[x_-(T)^\top\Sigma_k x_-(T)\right]& \leq  -\mathbb E\int_0^T \hspace{-0.25cm} x_-(t)^\top \bar C_k^\top \bar C_k x_-(t) dt \\   \nonumber
 &+\sigma_k  \mathbb E\int_0^T 2 x_{k, 2}(t) c_0(t) + \sum_{i, j=1}^q \left(2 \smat n_{i, 21} & n_{i, 22}\srix  x_{-}(t) +  c_i(t)\right)c_j(t) k_{ij} dt.
                                           \end{align}
Using the definitions of $c_i$ and $x_-$, we find an upper bound by replacing $x_-$ by $x_k$ in the last term, i.e.,
\begin{align}\label{important_est}
  \sum_{i, j=1}^q \left(2 \smat n_{i, 21} & n_{i, 22}\srix  x_-(t) +  c_i(t)\right)c_j(t) k_{ij} \leq \sum_{i, j=1}^q \left(2 \smat n_{i, 21} & n_{i, 22}\srix  x_{k}(t) +  c_i(t)\right)c_j(t) k_{ij}                                                                                         \end{align}
exploiting that $\sum_{i, j=1}^q c_i(t)c_j(t) k_{ij} \geq 0$ because $K=(k_{ij})$ is positive semidefinite.  
Secondly, we see that $\bar C_k x_- = \bar C_k x_k - \bar C_{k-1} x_{k-1}= \bar y_k-\bar y_{k-1}$ since $\bar C_k = \begin{bmatrix} \bar C_{k-1} &\star                                                                                                                                                                                                                                                            \end{bmatrix}$. Now, inserting these estimates  
into \eqref{estimate_used_again} implies
\begin{align}\label{first_result}
  \left\|\bar y_{k}-\bar y_{k-1}\right\|_{L^2_T}^2 \leq    
 \sigma_k  \mathbb E\int_0^T \hspace{-0.1cm} 2 x_{k, 2}(t) c_0(t) + \sum_{i, j=1}^q \left(2 \smat n_{i, 21} & n_{i, 22}\srix  x_{k}(t) +  c_i(t)\right)c_j(t) k_{ij} dt.
                                           \end{align}
We apply Lemma \ref{lemdgl} to  \eqref{plus} with $b(t) = 2 B v(t)$, $D= L_k$ and directly make use of \eqref{aux_reach_gram} providing\begin{align}\nonumber
&\mathbb E\left[x_+(T)^\top L_k x_+(T)\right] 
\leq \mathbb E\int_0^T - x_+(t)^\top L_k B_k B_k^\top L_k x_+(t) + 4 x_+(t)^\top L_k B_k v(t) dt\\ \label{estimate_used_again2}
& \quad-(\sigma_k + \sigma_k^{-1}) \mathbb E\int_0^T 2 x_{k, 2}(t) c_0(t)+ \sum_{i, j=1}^q \left(2 \smat n_{i, 21} & n_{i, 22}\srix  x_{+}(t) -  c_i(t)\right)c_j(t) k_{ij}   dt.
\end{align}
We observe that 
\begin{align}\label{estimate_on_xplus}
 \left(2 \smat n_{i, 21} & n_{i, 22}\srix  x_+(t) -  c_i(t)\right)c_j(t) = \left(2 \smat n_{i, 21} & n_{i, 22}\srix  x_{k}(t) +  c_i(t)\right)c_j(t)                                                                                      \end{align}
based on the definitions of $c_i$ and $x_+$ and furthermore find 
\begin{align}\label{controlest}
4\left\| v(t)\right\|_2^2 &\geq \left\|2 v(t)\right\|_2^2 -\left\|B_k^\top L_k x_+(t)- 2v(t)\right\|_2^2 \\ \nonumber
&=-x_+(t)^\top L_k B_kB_k^\top L_k x_+(t)+4 x_+(t)^\top L_k B_k v(t).
                \end{align}
Combining \eqref{estimate_used_again2} with \eqref{estimate_on_xplus} and \eqref{controlest} leads to 
\begin{align*}
 \mathbb E\int_0^T 2 x_{k, 2}(t) c_0(t)+ \sum_{i, j=1}^q \left(2 \smat n_{i, 21} & n_{i, 22}\srix  x_{k}(t) +  c_i(t)\right)c_j(t) k_{ij}   dt\leq \frac{4}{\sigma_k + \sigma_k^{-1}} \left\| v\right\|_{L^2_T}^2.
\end{align*}
This together with \eqref{first_result} gives us \begin{align*}
 \left\|\bar y_{k}-\bar y_{k-1}\right\|_{L^2_T}^2 \leq    
 4 \frac{\sigma_k}{\sigma_k + \sigma_k^{-1}} \left\| v\right\|_{L^2_T}^2 =
 4 \frac{\sigma_k^2}{\sigma_k^2 + 1} \left\| v\right\|_{L^2_T}^2.
\end{align*}
Inserting this into \eqref{remove_HSVs_individually}, it follows that $
 \left\|\begin{bmatrix} u-u_r\\ y-y_r\end{bmatrix}\right\|_{L^2_T} \leq 2 \sum_{k=r+1}^n \frac{\sigma_k}{\sqrt{\sigma_k^2 + 1}} \left\| v\right\|_{L^2_T}.
 $
It remains to calculate $\left\| v\right\|_{L^2_T}$ with $v(t)= B^\top Q x(t) + u(t)$. Based on \eqref{estimate_on_v}, we obtain
\begin{align*}
\mathbb E \left[x(T)^\top Q x(T)\right] 
=\int_0^T \mathbb E\left[-\left\|y(t)\right\|_2^2-\left\|u(t)\right\|_2^2 + \left\| B^\top Q x(t) + u(t)\right\|_2^2\right] dt,
\end{align*}
which provides the first claim of this theorem. If $u, x\in L^2$, then the limit as $T\rightarrow \infty$ of the above right hand side exists. Therefore, $\lim_{T\rightarrow\infty} \mathbb E \left[x(T)^\top Q x(T)\right]$  exists. Hence it is zero, otherwise it contradicts $x\in L^2$. Now, taking the limit as $T\rightarrow \infty$ in \eqref{bound_finite_time} yields the second claim.
\end{proof}

\section{Proof of Theorem \ref{main_error_bound2}}\label{appC}

\begin{proof}[Proof of Theorem \ref{main_error_bound2}]
Showing this result is more complex than the proof given in Appendix \ref{appB}. However, some basic steps are identical such that a similar notation will be used below. As before, let $(A_n, B_n, C_n, N_{i, n})$ be the balanced realization of \eqref{full_system}, i.e., the associated Gramians are identical and equal to $\Sigma_n$. Again, $A_{k}, B_{k}, C_{k}, N_{i, k}$ for $k= r, \dots, n-1$ are the respective submatrices of the balanced realization. They define the reduced system of dimension $k$ given by  
\begin{equation}\label{red_system3}
\begin{aligned}
 dx_k(t) &= [A_k x_k(t) +  {B_k} u(t)]dt + \sum_{i=1}^{q} {N}_{i, k} x_k(t) dW_i(t),\\
 y_k(t) &= C_k x_k(t),\quad t\geq 0.
\end{aligned}
\end{equation}
Setting $k=r$ now yields the reduced system \eqref{red_system} with $u_r= u$. Given that $k=n$, we obtain the balanced realization of \eqref{full_system} and hence $y_n=y$.
The inequality of Theorem
\ref{main_error_bound2} involves a scaled $L^2$-norm for which we can apply triangle inequality leading to
 \begin{align}\label{triangle_used}
\left( \mathbb E\int_0^T \expn^{-\beta t}\left\|y(t)-y_r(t)\right\|_{2}^2 dt \right)^{\frac{1}{2}}
 \leq \sum_{k=r+1}^n  \left( \mathbb E\int_0^T \expn^{-\beta t}\left\|y_k(t)-y_{k-1}(t)\right\|_{2}^2 dt \right)^{\frac{1}{2}}.                                                                                                                            \end{align}
 In order to proceed further, the error between $y_k$ and $y_{k-1}$ is analyzed. The associated matrix inequalities are derived from the ones for the balanced realization which are obtained by replacing $(A, B, C, N_i, P, Q)$ by $(A_n, B_n, C_n, N_{i, n}, \Sigma_n, \Sigma_n)$ in \eqref{obs_gram} and \eqref{reach_gram}. Evaluating the left upper $k\times k$ blocks of these (in)equalities yields \begin{subequations}\label{riccati_eq_bal}
\begin{align}\label{reach_gram_bal}
 A_k^\top \Sigma_k^{-1} + \Sigma_k^{-1} A_k + \sum_{i, j=1}^q N_{i, k}^\top \Sigma_k^{-1} N_{j, k} k_{ij} -C_k^\top C_k +\Sigma_k^{-1}B_k B_k^\top \Sigma_k^{-1}\leq 0,\\  \label{obs_gram_bal}
 A_k^\top \Sigma_k + \Sigma_k A_k + \sum_{i, j=1}^q N_{i, k}^\top \Sigma_k N_{j, k} k_{ij}+ C_k^\top C_k -\Sigma_k B_k B_k^\top \Sigma_k \leq 0
\end{align}
\end{subequations}
for $k=r, \dots, n$. We partition $x_k$, $N_{i, k}$ and $B_k$ like in \eqref{partition2}
and set $A_{k}= \begin{bmatrix} {A}_{k-1}&\star\\ 
{a}_{21}&\star\end{bmatrix}$. Below, the variables $x_-$ and $x_+$ are defined analogously to \eqref{xplusminus}. Based on \eqref{red_system3}, we find the respective equations by \begin{subequations}\label{plus_minus_sys2}
 \begin{align}\label{minus2}
 dx_-(t) &= [A_{k} x_-(t) +  \smat 0 \\ c_0(t)\srix]dt + \sum_{i=1}^{q}[ N_{i, k} x_-(t)+ \smat 0 \\ c_i(t)\srix] dW_i(t),\\ \label{plus2}
 dx_+(t) &= [A_{k} x_+(t) +  2 B_k u(t)-\smat 0 \\ c_0(t)\srix]dt + \sum_{i=1}^{q}[ N_{i, k} x_+(t)-\smat 0 \\ c_i(t)\srix]  dW_i(t),
\end{align}
\end{subequations}
where $c_0(t):=a_{21}x_{k-1}(t)+b_2 u(t)$ and $c_i(t):= n_{i, 21} x_{k-1}(t)$.
We apply Lemma \ref{lemdgl} to $\mathbb E\left[x_-(t)^\top\Sigma_k x_-(t)\right]$, $t\in[0, T]$, using \eqref{minus2} and exploit \eqref{obs_gram_bal} giving us
\begin{align}\label{estimate_used_again_again}
  \mathbb E\left[x_-(t)^\top\Sigma_k x_-(t)\right]& \leq \mathbb E\int_0^t \hspace{-0.25cm} x_-(s)^\top \Sigma_k B_k B_k^\top \Sigma_k x_-(s) ds   -\mathbb E\int_0^t \hspace{-0.25cm} x_-(s)^\top C_k^\top C_k x_-(s)ds \\   \nonumber
 &+\sigma_k  \mathbb E\int_0^t 2 x_{k, 2}(s) c_0(s) + \sum_{i, j=1}^q \left(2 \smat n_{i, 21} & n_{i, 22}\srix  x_{-}(s) +  c_i(s)\right)c_j(s) k_{ij} ds.                                        \end{align}
 We obtain that \begin{align*}                                x_-(s)^\top \Sigma_k B_k B_k^\top \Sigma_k x_-(s)=\left\|B_k^\top \Sigma_k^{\frac{1}{2}}\Sigma_k^{\frac{1}{2}} x_-(s)\right\|_{2}^2
  \leq    \left\|B_k^\top \Sigma_k^{\frac{1}{2}}\right\|_{2}^2 x_-(s)^\top \Sigma_k  x_-(s).                                                                                                                                                   \end{align*}
Since $B_k^\top \Sigma_k^{\frac{1}{2}}= B_n^\top \Sigma_n^{\frac{1}{2}}\begin{bmatrix}
 I_k\\
0                                                                 \end{bmatrix}
$, where $I_k$ is a $k\times k$ identity matrix, we have $\left\|B_k^\top \Sigma_k^{\frac{1}{2}}\right\|_{2}^2\leq \left\|B_n^\top \Sigma_n^{\frac{1}{2}}\right\|_{2}^2 = \left\|B^\top Q^{\frac{1}{2}}\right\|_{2}^2\leq \beta$. Therefore, we have
\begin{align*}                            x_-(s)^\top \Sigma_k B_k B_k^\top \Sigma_k x_-(s)
  \leq   \beta x_-(s)^\top \Sigma_k  x_-(s).                                                                                                                                                   \end{align*}
Moreover, we define $\alpha_k(t)= \mathbb E\int_0^t 2 x_{k, 2}(s) c_0(s) + \sum_{i, j=1}^q \left(2 \smat n_{i, 21} & n_{i, 22}\srix  x_{k}(s) +  c_i(s)\right)c_j(s) k_{ij} ds$ and see that $\alpha_k$ is an upper bound for the last integral in \eqref{estimate_used_again_again} taking \eqref{important_est} into account. We further observe that $C_k x_- = y_k - y_{k-1}$ such that 
 \eqref{estimate_used_again_again} becomes 
 \begin{align*}
  \mathbb E\left[x_-(t)^\top\Sigma_k x_-(t)\right]& \leq   \sigma_k\alpha_k(t)-\left\|y_k - y_{k-1}\right\|^2_{L^2_t}+\beta \int_0^t \hspace{-0.25cm} \mathbb E \left[x_-(s)^\top \Sigma_k x_-(s)\right] ds.                                        \end{align*}
We apply Lemma \ref{gronwall} resulting in \begin{align*}
  \mathbb E\left[x_-(t)^\top\Sigma_k x_-(t)\right] \leq   \sigma_k\alpha_k(t)-\left\|y_k - y_{k-1}\right\|^2_{L^2_t}+\beta \int_0^t \hspace{-0.15cm} [\sigma_k\alpha_k(s)-\left\|y_k - y_{k-1}\right\|^2_{L^2_s}]\expn^{\beta(t-s)}ds.                                        \end{align*}
Using integration by parts, we obtain \begin{align*}
 &                                  \beta \int_0^t \hspace{-0.25cm} [\sigma_k\alpha_k(s)-\left\|y_k - y_{k-1}\right\|^2_{L^2_s}]\expn^{\beta(t-s)}ds \\&= 
 \left[ -(\sigma_k\alpha_k(s)-\left\|y_k - y_{k-1}\right\|^2_{L^2_s} ) \expn^{\beta(t-s)}          \right]_{s=0}^t + \int_0^t \hspace{-0.15cm} [\sigma_k\dot\alpha_k(s)-\mathbb E\left\|y_k(s) - y_{k-1}(s)\right\|^2_{2}]\expn^{\beta(t-s)}ds.                       \end{align*}
Therefore, we have \begin{align*}
  \mathbb E\left[x_-(t)^\top\Sigma_k x_-(t)\right] \leq \int_0^t \hspace{-0.15cm} [\sigma_k\dot\alpha_k(s)-\mathbb E\left\|y_k(s) - y_{k-1}(s)\right\|^2_{2}]\expn^{\beta(t-s)}ds 
    \end{align*}
and hence, by setting $t=T$, we obtain \begin{align}\label{firstest_to_prove}
\mathbb E\int_0^T \hspace{-0.15cm} \left\|y_k(s) - y_{k-1}(s)\right\|^2_{2}\expn^{-\beta s}ds    \leq \sigma_k \int_0^T \hspace{-0.15cm} \dot\alpha_k(s)          \expn^{-\beta s}ds.
    \end{align}
In the following, an upper bound of the above right-hand side is found that depends on the control $u$. For that reason, we exploit \eqref{reach_gram_bal} after applying Lemma \ref{lemdgl} to find an expression for $ \mathbb E\left[x_+(t)^\top\Sigma_k^{-1} x_+(t)\right]$ based on \eqref{plus2}. Consequently, \begin{align}\nonumber
&\mathbb E\left[x_+(t)^\top \Sigma_k^{-1} x_+(t)\right] 
\leq \mathbb E\int_0^t x_+(s)^\top  C_k^\top C_k x_+(s) ds \\ \nonumber 
&\quad+\mathbb E\int_0^t - x_+(s)^\top \Sigma_k^{-1} B_k B_k^\top \Sigma_k^{-1} x_+(s) + 4 x_+(s)^\top \Sigma_k^{-1} B_k u(s) ds\\ \label{estimate_used_again2_again2}
& \quad- \sigma_k^{-1} \mathbb E\int_0^t 2 x_{k, 2}(s) c_0(s)+ \sum_{i, j=1}^q \left(2 \smat n_{i, 21} & n_{i, 22}\srix  x_{+}(s) -  c_i(s)\right)c_j(s) k_{ij}   ds.
\end{align}
With the same argument like in \eqref{controlest}, it can be shown that\begin{align*}
4\left\| u(s)\right\|_2^2 &\geq -x_+(s)^\top \Sigma_k^{-1} B_kB_k^\top \Sigma_k^{-1} x_+(s)+4 x_+(s)^\top \Sigma_k^{-1} B_k u(s).
    \end{align*}
Using the definitions of $x_+$ and $c_i$, it immediately follows that $\mathbb E\int_0^t 2 x_{k, 2}(s) c_0(s)+ \sum_{i, j=1}^q \left(2 \smat n_{i, 21} & n_{i, 22}\srix  x_{+}(s) -  c_i(s)\right)c_j(s) k_{ij}   ds = \alpha_k(t)$. Inserting these insights into \eqref{estimate_used_again2_again2}, we obtain \begin{align*}
\mathbb E\left[x_+(t)^\top \Sigma_k^{-1} x_+(t)\right] 
\leq \mathbb E\int_0^t x_+(s)^\top  C_k^\top C_k x_+(s) ds + 4 \left\| u\right\|_{L^2_t}^2-\sigma_k^{-1} \alpha_k(t).
\end{align*}
Since it holds that
$\left\|C_k \Sigma_k^{\frac{1}{2}}\right\|_{2}^2\leq \left\|C_n \Sigma_n^{\frac{1}{2}}\right\|_{2}^2 = \left\|C P^{\frac{1}{2}}\right\|_{2}^2\leq \beta$, we have \begin{align*}
\mathbb E\left[x_+(t)^\top \Sigma_k^{-1} x_+(t)\right] 
\leq \beta \mathbb E\int_0^t x_+(s)^\top  \Sigma_k^{-1} x_+(s) ds + 4 \left\| u\right\|_{L^2_t}^2-\sigma_k^{-1} \alpha_k(t).
\end{align*}
Lemma \ref{gronwall} now delivers\begin{align*}
\mathbb E\left[x_+(t)^\top \Sigma_k^{-1} x_+(t)\right] 
\leq \beta \int_0^t[4 \left\| u\right\|_{L^2_s}^2-\sigma_k^{-1} \alpha_k(s)] \expn^{\beta(t-s)} ds + 4 \left\| u\right\|_{L^2_t}^2-\sigma_k^{-1} \alpha_k(t).
\end{align*}
Again, integration by parts leads to \begin{align*}
\mathbb E\left[x_+(t)^\top \Sigma_k^{-1} x_+(t)\right] 
\leq \int_0^t[4 \mathbb E\left\| u(s)\right\|_{2}^2-\sigma_k^{-1} \dot \alpha_k(s)] \expn^{\beta(t-s)} ds.
\end{align*}
Consequently, we find that 
\begin{align*}
\int_0^T \dot \alpha_k(s) \expn^{ -\beta s} ds\leq 4\sigma_k \mathbb E  \int_0^T\left\| u(s)\right\|_{2}^2 \expn^{-\beta s} ds.
\end{align*}
Using this estimate for \eqref{firstest_to_prove}, the result follows from \eqref{triangle_used}.
\end{proof}

\section*{Acknowledgments}
 MR is supported by the DFG via the individual grant ``Low-order approximations for large-scale problems arising in the context of high-dimensional
PDEs and spatially discretized SPDEs''-- project number 499366908.

\bibliographystyle{plain}
\bibliography{refs}

\end{document}